# WHY DIFFEOLOGY?

PATRICK IGLESIAS-ZEMMOUR

Dedicated to the memory of Vladimir Igorevich Arnold

ABSTRACT. Diffeology extends differential geometry to spaces beyond smooth manifolds. This paper explores diffeology's key features and illustrates its utility with examples including singular and quotient spaces, and applications in symplectic geometry and prequantization. Diffeology provides a natural and effective framework for handling complexities from singularities and infinite-dimensional settings.

## INTRODUCTION: THE NEED FOR A BROADER FRAMEWORK

Why bother with diffeology? Classical differential geometry works well enough, doesn't it? The answer, to put it bluntly, is no. From the awkwardness of handling singularities to the limitations imposed by infinite-dimensional spaces, the traditional framework is riddled with inadequacies. Diffeology offers not just a patch, but a fundamental overhaul, a comprehensive theory, providing a more natural and powerful foundation for smoothness. This paper will convince you that it's a theory worth learning.

Classical differential geometry often falls short when dealing with quotient spaces that are not manifolds (such as those arising from non-free group actions or dynamical systems like the irrational torus $T_\alpha$), infinite-dimensional spaces (like spaces of functions or groups of diffeomorphisms), and singularities. To move beyond ad hoc modifications and heuristic constructions, a more general framework is required. Diffeology provides an effective and minimal response to this need by abandoning the compulsory Euclidean local model. Smoothness is instead defined intrinsically by a family of parametrizations (maps from Euclidean domains) called 'plots', subject to a few natural axioms.

This approach allows us to treat such singular spaces as *autonomous geometrical objects*, possessing their own inherent smooth structure, independent of any embedding or originating space. Following Felix Klein's Erlangen Program [Kle72], a geometry is fundamentally characterized by a space and a group of transformations acting upon it,

---

*Date*: November 29, 2025.

2020 *Mathematics Subject Classification.* Primary 58A40; Secondary 57S15, 53D50, 53D20, 55R10, 58A05.

I gratefully acknowledge the Einstein Institute of the Hebrew University of Jerusalem for their continued support, and the developers of the AI assistant Gemini (Google) for their invaluable assistance in writing this paper.

focusing on the properties invariant under this group. Diffeology fits this paradigm perfectly: every diffeological space X comes equipped with its group of diffeomorphisms, Diff(X), [1] and this is my definition of *Differential Geometry*.The study of X's geometry, in this Kleinian sense, involves analyzing this group and its associated invariants. The non-triviality of these invariants, even for spaces with trivial topology like $T_\alpha$, reveals a rich underlying geometric structure.

As a category, diffeology exhibits remarkable properties: it is complete, co-complete, and Cartesian closed. This means it is stable under all set-theoretic operations, and crucially, that sets of smooth maps ($\mathscr{C}^\infty(X, X')$) are themselves objects within the category. These categorical properties enable simple constructions, but their true power lies in the theory's ability to provide a fine analysis of singularities and singular spaces, exemplified by the irrational torus [DIZ83], whose non-trivial diffeological invariants encode arithmetic properties of its defining slope [IZL90] [PIZ25a]. It is the conjunction of these categorical properties and the ability to capture fine geometric structure that distinguishes diffeology.

The mathematical framework of diffeology, which provides a powerful tool for studying such spaces, was originally laid out by Jean-Marie Souriau in his seminal work "Groupes Différentiels" [Sou80]. Souriau's axiomatic approach to defining smoothness of infinite dimensional groups through plots paved the way for a consistent and flexible theory applicable to a wide range of spaces beyond traditional manifolds. It is worth noting that the seeds of this idea were also present in the work of Kuo-Tsai Chen "Iterated Path Integrals" [Che77], which explored similar concepts with a slightly different kind of plots and a different goal.

Of course, diffeology is not the only candidate for generalization; there are also noncommutative geometry, Frölicher spaces, Sikorski differentiable spaces, derived differential geometry, and synthetic differential geometry. However, Sikorski and Frölicher spaces render the irrational torus trivial,[2] while Frölicher spaces embed faithfully into {Diffeology} as reflexive spaces. Recent work has exhibited a bridge between diffeology and noncommutative geometry, a connection that deserves further exploration. The relationships with derived differential geometry or synthetic differential geometry are, to my knowledge, yet to be fully clarified.

This paper serves as an invitation to diffeology. It focuses on presenting the core concepts and applications and does not delve into technical constructions or proofs, except in a couple of cases that deserve special attention. For a deeper understanding, interested readers are encouraged to consult the textbook on diffeology [PIZ13]. For a broader range of examples, please consult [PIZ25]. This is not an inventory of all that has been achieved in diffeology, nor a complete list of contributions and contributors;

---

[1] While the global group Diff(X) captures the overall symmetries of the space, singularities are a fundamentally local phenomenon. To properly classify these local geometric structures, the Kleinian paradigm is naturally extended to consider the action of the richer *pseudogroup of local diffeomorphisms*. As we will see, the orbits of this pseudogroup — the Klein strata — provide the definitive classification of a space's intrinsic singularities.

[2] In the sense that its structure is degenerate.

rather, it highlights its role in creating new avenues to address limitations in differential geometry, with the expectation of continued development and broader application.

This paper will not only outline the foundational principles of diffeology but will also demonstrate its utility through some decisive applications. In particular, we will show how this framework provides the definitive dictionary for the intrinsic geometry of singular orbit spaces [GIZ25], how it resolves foundational paradoxes in geometric quantization through a novel path-space formulation [PIZ25b]. These examples serve as irrefutable evidence that diffeology is not merely a generalization, but a necessary tool for a deeper understanding of the smooth world.

A last word of introduction before starting the hard work: I asked Gemini Google about diffeology, this is what it replies to me.

> Is 'smoothness' a property of spaces, or a construction of our mathematical tools? Classical differential geometry often blurs this line, leading to artificial constraints. Diffeology dares to redefine smoothness from the ground up, asking: what is the most natural way to describe how we probe a space with smooth maps? The result is a surprisingly powerful and flexible theory that liberates us from the limitations of manifolds. So, why would you need diffeology? Perhaps because it offers a more fundamental and less biased view of the smooth world.

I would argue that if an AI, trained on the vast corpus of mathematical knowledge, identifies the limitations of classical differential geometry and articulates the need for a more fundamental approach like diffeology, then it speaks to the inherent validity and importance of this perspective.

## 1. THE CATEGORY {DIFFEOLOGY}: FOUNDATIONAL PROPERTIES

Everything starts with *parametrizations* which are the primary tool for investigating diffeological spaces. A parametrization of a set X is any map P from a Euclidean domain (an open subset of a Euclidean space) into X. A diffeology on X is a choice of parametrizations considered *smooth*. These smooth parametrizations are called *plots*.

**1. What is a diffeology?** A diffeology on a set X is a collection $\mathcal{D}$ of plots satisfying the following axioms, considered the minimal conditions for a meaningful theory of smoothness:

(1) *Covering*. Every constant parametrization $r \mapsto x$ from any Euclidean domain into X is a plot. This ensures that every point $x \in \mathrm{X}$ can be "reached smoothly."
(2) *Smooth Compatibility*. Composition with smooth maps from Euclidean domains preserves smoothness: If $\mathrm{P} : \mathrm{U} \to \mathrm{X}$ is a plot and $\mathrm{F} : \mathrm{V} \to \mathrm{U}$ is a smooth map, then $\mathrm{P} \circ \mathrm{F} : \mathrm{V} \to \mathrm{X}$ is a plot.
(3) *Locality*. Smoothness is a local property: Let $\mathrm{P} : \mathrm{U} \to \mathrm{X}$ be a parametrization such that for every $r \in \mathrm{U}$, there is an open neighborhood $\mathrm{V} \subset \mathrm{U}$ of $r$ where $\mathrm{P} \restriction \mathrm{V}$ is a plot. Then P is a plot.

A *diffeological space* is a set equipped with a diffeology. It is crucial to note that the underlying set initially possesses no additional structure; a topology may be subsequently derived, but it is a secondary construction and can even be trivial, as in the case of the irrational torus.

Euclidean domains themselves are the first and most elementary examples of diffeological spaces, equipped with the standard $\mathscr{C}^\infty$ diffeology where plots are ordinary smooth maps.

Diffeologies are partially ordered by inclusion, called *fineness*. On any set, the finest diffeology is the *discrete diffeology* (plots are locally constant parametrizations), and the coarsest is the *coarse* or *vague diffeology* (all parametrizations are plots). The set of diffeologies on a set forms a complete lattice, meaning every subset of diffeologies has an infimum and a supremum. This property is useful for discriminate diffeologies satisfying a given property.

**2. Smooth maps.** A map $f : \mathrm{X} \to \mathrm{X}'$ between diffeological spaces is *smooth* if it maps plots to plots by composition: for every plot P in X, $f \circ \mathrm{P}$ is a plot in X′. The set of smooth maps from X to X′ is denoted by $\mathscr{C}^\infty(\mathrm{X}, \mathrm{X}')$.

The plots in a diffeological space are exactly its smooth parametrizations.

By associativity, the composition of smooth maps is again a smooth map, which ensures that the composition of smooth maps is well-defined, and that the diffeological spaces form a category. We denote it by {Diffeology}.

The isomorphisms of this category are any smooth bijection $f$, such that its inverse $f^{-1}$ is also smooth.

**3. Example: The Irrational Torus.** The irrational torus was one of the first non-trivial diffeological spaces studied [DIZ83]. It remains a fundamental example illustrating the capabilities of diffeology and was a key motivation for its development. The power of diffeology lies in the combination of its categorical properties (complete, co-complete, and Cartesian closed) with the highly non-trivial nature of the irrational torus. This non-triviality allows for generalizing the prequantum bundle construction for any parasymplectic form[3] on manifolds with any group of periods [PIZ95]. The fiber of the generalized prequantum bundle, the *torus of periods*, is the quotient of the real line by the *group of periods*. When the periods have more than one generator, it is a generic irrational torus, not a manifold, and consequently, the prequantum bundle is also not a manifold. This generalization reframes the traditional prequantum bundle (where the group of periods has a single generator) as a special case within an infinite spectrum of possibilities, a spectrum that is absent in smooth categories where irrational tori are trivial.

For clarity, we consider the torus T in **C** as the subset of complex numbers with modulus 1. A plot in T is any smooth map $r \mapsto z(r) = x(r) + i\,y(r)$ in **C** such that $|z(r)| = 1$ for all

[3]We use the term 'parasymplectic' to refer to a general closed 2-form on a diffeological space, without assuming any non-degeneracy condition. This distinguishes it from the stronger notions of presymplectic or symplectic forms.

$r$. A plot in $\mathrm{T}^2$ is a pair of plots in T, $r \mapsto (z_1(r), z_2(r))$, defining the 2-torus as a standard diffeological space.

Let $\alpha \in \mathbf{R} - \mathbf{Q}$ be an irrational number, and consider the free action of the additive group $\mathbf{R}$ on $\mathrm{T}^2$ defined by:

$$\underline{t}(z_1, z_2) = \left(e^{2i\pi t} z_1, e^{2i\pi\alpha t} z_2\right)$$

The irrational torus $\mathrm{T}_\alpha$ is defined as the quotient of $\mathrm{T}^2$ by this action of $\mathbf{R}$ [DIZ83]. Since $\mathrm{T}^2$ is an Abelian group and the subset

$$\mathscr{S}_\alpha = \left\{\left(e^{2i\pi t}, e^{2i\pi\alpha t}\right) \mid t \in \mathbf{R}\right\}$$

is a subgroup, the irrational torus is the group quotient

$$\mathrm{T}_\alpha = \mathrm{T}^2 / \mathscr{S}_\alpha.$$

Let class : $\mathrm{T}^2 \to \mathrm{T}_\alpha$ be the projection. The diffeology of $\mathrm{T}_\alpha$ is defined by the parametrizations $r \mapsto \tau(r)$ such that, for every point $r$ in its domain, there exists a neighborhood V of $r$ and a plot $r \mapsto (z_1(r), z_2(r))$ in $\mathrm{T}^2$ such that

$$\text{For all } r \in \mathrm{V}, \quad \tau(r) = \mathrm{class}(z_1(r), z_2(r)).$$

This is illustrated by the diagram of a plot P:

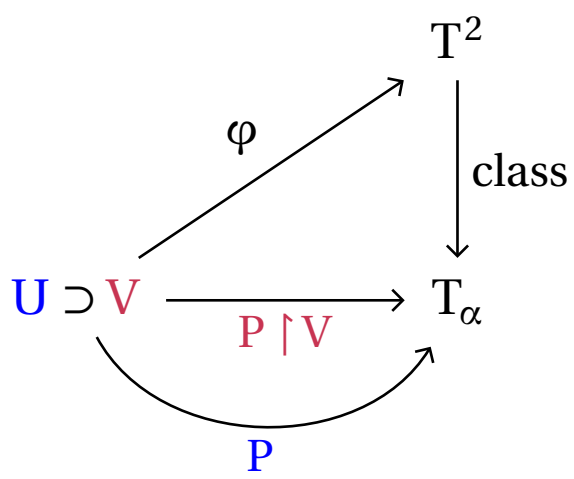

The irrational torus, the quotient of an Abelian Lie group by a subgroup, is a *diffeological group* but not a Lie group, because $\mathrm{T}_\alpha$ inherits the trivial topology and is not a manifold.

The irrational torus satisfies the following properties (op. cit.):

Fact 1 The irrational torus $\mathrm{T}_\alpha$ is diffeomorphic to the quotient $\mathbf{R}/(\mathbf{Z} + \alpha\mathbf{Z})$.

Fact 2 The projection $\pi : \mathbf{R} \to \mathrm{T}_\alpha \simeq \mathbf{R}/(\mathbf{Z} + \alpha\mathbf{Z})$ is its universal covering, unique up to an isomorphism.

Fact 3 The first homotopy group of $\mathrm{T}_\alpha$ is isomorphic to $\mathbf{Z} + \alpha\mathbf{Z} \subset \mathbf{R}$.

Fact 4 Two irrational tori $\mathrm{T}_\alpha$ and $\mathrm{T}_\beta$ are diffeomorphic if they are conjugate by a unimodular transformation; that is, if there exists

$$\begin{pmatrix} a & b \\ c & d \end{pmatrix} \in \mathrm{GL}(2, \mathbf{Z}) \text{ such that } \beta = \frac{a\alpha + b}{c\alpha + d}.$$

Fact 5. Every diffeomorphism of $\mathrm{T}_\alpha$ is the projection of some affine map $x \mapsto \lambda x + \mu$ where $\mu$ is any number, and $\lambda$ belongs to a subgroup of the multiplicative group of non-zero real numbers, isomorphic to the group of connected components of the group of diffeomorphisms $\mathrm{Diff}(\mathrm{T}_\alpha)$. And we have:[4]

$$\pi_0(\mathrm{Diff}(\mathrm{T}_\alpha)) \simeq \begin{cases} \{\pm 1\} \times \mathbf{Z}, & \text{if } \alpha \text{ is a quadratic number;} \\ \{\pm 1\} & \text{otherwise.} \end{cases}$$

Thus, the group $\mathrm{Diff}(\mathrm{T}_\alpha)$ is either isomorphic to $\{\pm 1\} \times \mathbf{Z} \times \mathrm{T}_\alpha$ or $\{\pm 1\} \times \mathrm{T}_\alpha$, depending on whether $\alpha$ is quadratic or not. The identity component $\mathrm{Diff}(\mathrm{T}_\alpha)^0$ is isomorphic to $\mathrm{T}_\alpha$ itself, acting by multiplication.

There is a generalization of this result when $\alpha = (\alpha_i)_{i=1}^n$ is a linear 1-form on $\mathbf{R}^n$ and the hyperplane $\ker(\alpha)$ is irrational [IZL90]. In this case, $\pi_0(\mathrm{Diff}(\mathrm{T}_\alpha))$ is the group of units of an order of a finite field extension $\mathrm{K}_\alpha$ of $\mathbf{Q}$ associated with the coefficient $\alpha_i$.

The properties listed above for the irrational torus demonstrate how diffeology can be used to analyze such spaces.

**4. Example: Manifolds, the role of the D-topology.** Diffeology redefines manifolds using local diffeomorphisms, which, along with the D-topology, are central to modeling diffeological spaces.

DEFINITION (LOCAL SMOOTH MAPS AND D-TOPOLOGY). *A map $f : \mathrm{X} \supset \mathrm{A} \to \mathrm{X}'$ from a subset* A *of a diffeological space* X *to a diffeological space* X′ *is* locally smooth *if for every plot* $\mathrm{P} : \mathrm{U} \to \mathrm{X}$ *in* X, $f \circ \mathrm{P}$ *is a plot in* X′. *This implies* $\mathrm{P}^{-1}(\mathrm{A}) \subset \mathrm{U}$ *is open (a Euclidean domain). The subsets* A *of* X *where* $\mathrm{P}^{-1}(\mathrm{A})$ *is open for all plots in* X *form a topology on* X, *called the* D-topology. *Smooth maps are D-continuous, and the D-topology is the finest topology making plots continuous.*

The trivial D-topology of the irrational torus illustrates diffeology's ability to capture smooth local structure even in such a critical case.

The composition of local smooth maps is locally smooth (including the empty map), and smooth maps define a category with objects the D-open subsets of diffeological spaces [IZL18]. The isomorphisms of this category are *local diffeomorphisms*.

DEFINITION (LOCAL DIFFEOMORPHISMS). *A local diffeomorphism from a diffeological space* X *to* X′ *is an injective local smooth map* $f : \mathrm{X} \supset \mathrm{A} \to \mathrm{X}'$ *with a local smooth inverse* $f^{-1} : \mathrm{X}' \supset f(\mathrm{A}) = \mathrm{A}' \to \mathrm{X}$.

Local diffeomorphisms are instrumental in defining subcategories within diffeology based on model spaces, such as manifolds:

DEFINITION (MANIFOLDS). *A manifold is a diffeological space locally diffeomorphic to a Euclidean space.*

Manifolds form a full and faithful subcategory of {Diffeology}, denoted {Manifolds}.

**5. Operations on Diffeology.** The diffeology of $\mathrm{T}_\alpha$ exemplifies a quotient diffeology, one of the four fundamental categorical operations on diffeologies: *sums, products, subspaces*, and *quotients*. These operations ensure the stability of the category and are fundamental for constructing new diffeological spaces from existing ones, expanding the theory's scope. The building blocks are Euclidean domains, and these operations are the tools upon which differential geometry is built. Using these operations ensures natural and consistent constructions, avoiding the ad-hoc modifications often required in the traditional approach, which lacks this stability.

[4]A quadratic number is a solution of a quadratic equation with integer coefficients.

1. *Sums.* The sum (or disjoint union) X of a family $\{X_i\}_{i\in\mathcal{I}}$ of diffeological spaces,
$$X = \coprod_{i\in\mathcal{I}} X_i = \{(i,x) \mid i \in \mathcal{I} \text{ and } x \in X_i\},$$
has the *diffeology sum*: the finest diffeology making each inclusion $j_i : X_i \to X$ smooth. Its plots are $r \mapsto (i(r), x(r))$ where $i(r)$ is locally constant and $r \mapsto x(r)$ is a plot in $X_i$.
2. *Products.* The product X of a family $\{X_i\}_{i\in\mathcal{I}}$ of diffeological spaces,
$$X = \prod_{i\in\mathcal{I}} X_i = \{(x_i)_{i\in\mathcal{I}} \mid x_i \in X_i\},$$
has the *diffeology product*: the coarsest diffeology making each projection $\pi_i : X \to X_i$ smooth. Its plots are $r \mapsto (x_i(r))$ where each $r \mapsto x_i(r)$ is a plot in $X_i$.
3. *Subspaces.* A subset $X'$ of a diffeological space X inherits the *subset diffeology*: plots in $X'$ are plots in X. This is the coarsest diffeology making the injection $j : X' \to X$ smooth. $X'$ equipped with this diffeology is a subspace.
4. *Quotients.* A surjective map $\pi : X \to Q$ from a diffeological space X to a set Q defines the *quotient diffeology* on Q: plots are $r \mapsto q(r)$ that lift locally smoothly in X. That is, for every $r$ in the domain, there is a neighborhood V and a plot $r \mapsto x(r)$ in X such that $q(r) = \pi(x(r))$ for all $r \in V$. This is the finest diffeology making $\pi$ smooth. The quotient $X/\sim$ of a diffeological space X by an equivalence relation $\sim$ inherits a natural quotient diffeology.[5]

In addition to these operations, the smooth maps between diffeological spaces form a diffeological space under the *functional diffeology*, making spaces of smooth maps themselves objects within diffeology.

5. *Functional Diffeology (Spaces of Smooth Maps).* A plot $r \mapsto f_r$ in $\mathcal{C}^\infty(X, X')$ (where X and $X'$ are diffeological spaces) is a parametrization on U such that $(r, x) \mapsto f_r(x)$, defined on $U \times X$, is smooth.

*Example 1 (The Constructive Principle).* Every diffeological space can be constructed from Euclidean domains by successively applying some of these operations. For instance, the irrational torus is built from four copies of **R** (with their standard diffeology, i.e., $\mathbf{C} \times \mathbf{C}$) by taking the subspace $T \times T \subset \mathbf{C} \times \mathbf{C}$, then quotienting $T^2$ by the orbits of **R** under the irrational solenoid $\mathcal{S}_\alpha$ to obtain $T_\alpha = T^2/\mathcal{S}_\alpha$. This exemplifies the constructive mode: product, subspace, quotient.

*Example 2 (The Universal Principle).* Every diffeological space X is the quotient of the sum $\mathcal{N}$ of the domains of any generating family[6] (the *Nebula* of the family) by the partition defined by the inverse images of the evaluation function $(F, r) \mapsto F(r)$ from $\mathcal{N}$ to X. Thus, every diffeological space can be constructed from Euclidean domains using sums and quotients. Actually, considering Euclidean domains as subspaces of products of $\mathbf{R}$, we use all four operations in a sequence: product, subspace, sum and quotient.

[5] The quotient diffeology depends only on the partition of X into classes, not on how the partition is obtained. For example, if two different group actions on X have the same orbits, the quotient diffeology is the same. This highlights how diffeologies handle quotients.

[6] A family of plots generating the diffeology: the finest diffeology containing these plots (the intersection of all diffeologies containing them).

This example reveals a profound and beautiful universality theorem at the heart of the theory. It makes the astonishing claim that *every* diffeological space, no matter how abstract, can be realized in the same way: as a simple quotient of a classical (though typically vast) smooth simple manifold: a sum of Euclidean domains. This is a radical unification. It provides a single, universal model for all of smoothness. It says that, at the most fundamental level, the only “source" of smoothness you need is $\mathbf{R}^n$. Everything else is just a shadow or a projection of that source. This is the ultimate answer to the question “Where do diffeological spaces come from?" They all come from one place. This is the theory's “creation myth”.

## 2. Navigating Singularities with Diffeology

**6. How to do differential geometry with a singular object.** Differential geometry essentially consists of comparing spaces and their properties. Returning to the irrational torus, let us illustrate how diffeology, an extended version of differential geometry, works in a concrete situation. Assume that $f : \mathrm{T}_\alpha \to \mathrm{T}_\beta$ is a diffeomorphism. Let $\pi_\alpha : \mathbf{R} \to \mathrm{T}_\alpha$ and $\pi_\beta : \mathbf{R} \to \mathrm{T}_\beta$ be the two projections from their universal coverings. Since $\pi_\alpha$ is a plot in $\mathrm{T}_\alpha$, it has a local lifting $\phi : \mathrm{I} \to \mathbf{R}$ over $f$, that is, $f \circ \pi_\alpha \restriction \mathrm{I} = \pi_\beta \circ \phi$.

$$\begin{array}{ccc} \mathbf{R} \supset \mathrm{I} & \xrightarrow{\ \phi\ } & \mathbf{R} \\ \pi_\alpha \big\downarrow & & \big\downarrow \pi_\beta \\ \mathrm{T}_\alpha & \xrightarrow[\ f\ ]{} & \mathrm{T}_\beta \end{array}$$

Hence, for all $n, m$ in $\mathbf{Z}$, there exist $n', m'$ in $\mathbf{Z}$ such that

$$\phi(x + n + m\alpha) = \phi(x) + n' + m'\beta.$$

By differentiating this equality, and thanks to the fact that $\mathbf{Z} + \alpha\mathbf{Z}$ is dense in $\mathbf{R}$, we deduce that the map $\phi$ is affine:

$$\phi(x) = \lambda x + \mu,$$

with $\lambda, \mu \in \mathbf{R}$. Substituting this expression into the equality above, we deduce that there exist four integers $a, b, c, d$, such that

$$\beta = \frac{a\alpha + b}{c\alpha + d}.$$

And since $f$ is bijective, we conclude that this correspondence is a unimodular transformation. And this is how we conduct an analysis in diffeology, in this case.

**7. Quotient Spaces and Dimension.** Let us consider another type of singularity, and examine the quotient space $\Delta_n = \mathbf{R}^n/\mathrm{O}(n)$, where $\mathrm{O}(n)$ is the orthogonal group. This space can be identified with the half-line $[0, \infty[$ equipped with the quotient diffeology induced by the standard diffeology of $\mathbf{R}^n$ and the map $\mathrm{sq} : \mathbf{R}^n \to [0, \infty[$, where $\mathrm{sq}(x) = \|x\|^2$. The plots in $[0, \infty[$ that constitute the diffeology of $\Delta_n$ are non-negative parametrizations $r \mapsto t(r)$ that can be expressed locally as $t(r) = \|x(r)\|^2$, where $r \mapsto x(r)$ is a smooth parametrization in $\mathbf{R}^n$.

The D-topology of $\Delta_n$ coincides with the standard topology of $[0,\infty[ \subset \mathbf{R}$ for all $n$. We now consider the diffeology.

THEOREM (DIFFEOLOGY OF THE HALF-LINES). *For any two distinct integers $n \neq m$, the half-lines $\Delta_n$ and $\Delta_m$ are not diffeomorphic. They are also not diffeomorphic to $\Delta_\infty = [0,\infty[$ equipped with the subset diffeology.*

The proof of this theorem relies on the concept of *dimension* in diffeology, and the *dimension map* [PIZ07]. For a chosen point $x \in \mathrm{X}$, the germ of the diffeology at $x$ is the set of all plots $\mathrm{P} : \mathrm{U} \to \mathrm{X}$ in X centered at $x$, i.e., such that $\mathrm{P}(0) = x$, where $0 \in \mathrm{U}$. The germ of a diffeology at a point $x$ admits the same notion of generating family as the diffeology of the whole space. A generating family $\mathscr{F}_x$ is a set of centered plots that generates the germ of the diffeology at $x$. The dimension of such a family $\mathscr{F}_x$ is defined by

$$\dim(\mathscr{F}_x) = \sup_{\mathrm{P} \in \mathscr{F}_x} \dim(\mathrm{dom}(\mathrm{P})).$$

The dimension of X at the point $x$ is then defined as:

$$\dim_x(\mathrm{X}) = \inf_{\mathscr{F}_x \in \mathrm{Gen}_x(\mathrm{X})} \dim(\mathscr{F}_x),$$

where $\mathrm{Gen}_x(\mathrm{X})$ is the set of all generating families of X at $x$.[7] Intuitively, $\dim_x(\mathrm{X})$ represents the smallest dimension of a Euclidean domain needed to locally generate the diffeology around $x$. The map $x \mapsto \dim_x(\mathrm{X})$ is called the dimension map. Since diffeomorphisms transform generating families into generating families, the dimension map is a diffeological invariant.

THEOREM (DIMENSION INVARIANT). *A diffeomorphism between two diffeological spaces maps each point of the first space to a point of the second space with the same dimension.*

Applying this to the half-lines, we can compute their dimension maps. For $\Delta_n$, the origin (corresponding to $t = 0$) has dimension $n$, reflecting the $n$ degrees of freedom in $\mathbf{R}^n$ collapsing at the origin under the quotient. Any other point $t > 0$ has dimension 1, since locally, the half-line behaves like a standard interval.

$$\dim_0(\Delta_n) = n \quad \text{and} \quad \dim_t(\Delta_n) = 1, \quad \text{for all } t > 0.$$

This holds for $n = \infty$ as well. Consequently, $\Delta_n \not\simeq \Delta_m$ if $n \neq m$, as their dimension maps are distinct. In this case, the dimension is sufficient to distinguish between non-diffeomorphic spaces, or points within a space.

The dimension in diffeology can be useful in the context of quotients of manifolds. For example, consider the quotient of a symplectic manifold $(\mathrm{M}^{2n}, \omega)$ by an effective action of the torus $\mathrm{T}^n$. The dimension map of the quotient space $\mathrm{Q}_n = \mathrm{M}^{2n}/\mathrm{T}^n$ is related to the depth, a concept often introduced ad hoc to describe the topological structure of

[7] This is what is usually called a "minimax" definition.

$Q_n$.[8] Specifically,

$$\dim_x(Q_n) = n + \text{depth}(x).$$

In the local model, the quotient $Q_n = \mathbf{C}^n/T^n$ indeed has the D-topology of a corner $[0,\infty[^n$. The depth corresponds to the number of coordinates that are zero. The dimension at $x$ is determined by this formula, avoiding the ad-hoc introduction of the depth as a heuristic in this construction, see [GIZ25] for a discussion on that question. We will see later that this is part of a natural stratification of diffeological spaces, the Klein stratification.

**8. Klein Stratification.** The analysis of the half-lines reveals that the singular origins of two half lines with different dimensions cannot be mapped by a diffeomorphism, but they can by a homeomorphism. The case of the square presents a similar scenario: consider the square $\text{Sq} = [0,1]^2 \subset \mathbf{R}^2$ equipped with the subset diffeology. Topologically, it is equivalent to a disc $D = \{x \in \mathbf{R}^2 \mid \|x\| \leq 1\}$, and a homeomorphism can map a point on the boundary to any other boundary point. However, diffeomorphisms behave differently. While a diffeomorphism of the disc can map any boundary point to any other boundary point, this is not the case for the square. The orbits of diffeomorphisms exhibit a different structure:

$$\begin{cases} D/\text{Diff}(D) = \{\text{boundary, interior}\}, \\ \text{Sq}/\text{Diff}(\text{Sq}) = \{\text{corners, edges, interior}\}. \end{cases}$$

Diffeologies can discern structures that are topologically equivalent, highlighting the finer resolution provided by diffeomorphisms. This aligns perfectly with the Kleinian perspective: the group Diff(Sq) reveals more geometric structure than the group of homeomorphisms. This motivates the following definition.

DEFINITION (KLEIN STRATA). *The* Klein strata *of a diffeological space are the orbits of its group of* local *diffeomorphisms.*

Several remarks are pertinent to this definition. First, Klein strata can be interpreted as inherent singularities of the space, revealed by the action of its local symmetries (local diffeomorphisms). Second, we use local diffeomorphisms rather than global ones because singularities are fundamentally local phenomena. Finally, the term "Klein strata" honors Felix Klein [Kle72], who pioneered the redefinition of geometry as the study of a space acted upon by a group of transformations and the invariants thereof. The Klein stratification provides a concrete realization of this program within diffeology, partitioning the space based on local symmetry types.

THEOREM (KLEIN STRATIFICATION). *The Klein strata of a diffeological space constitute an internal stratification, termed the Klein stratification, which satisfies the fundamental frontier condition: the D-topology closure of any stratum is a union of strata.*

Traditional stratification theory typically begins with a topological space and imposes additional requirements on the strata beyond the frontier condition. For instance, strata are often required to possess a manifold structure compatible with the topology, and to be locally closed to ensure the set of strata inherits a Hausdorff topology. These conditions are unnecessary in diffeology for two main reasons: first, strata naturally inherit the subset diffeology, remaining within the category of diffeological spaces. Second, as illustrated by the irrational torus, diffeology does not necessitate a Hausdorff topology for the set of strata; it can be trivial or partially trivial and can possess a highly non-trivial diffeology.

[8] In particular in the yet unpublished paper: *Classification of Locally Standard Torus Actions* by Yael Karshon and Shintarô Kuroki.

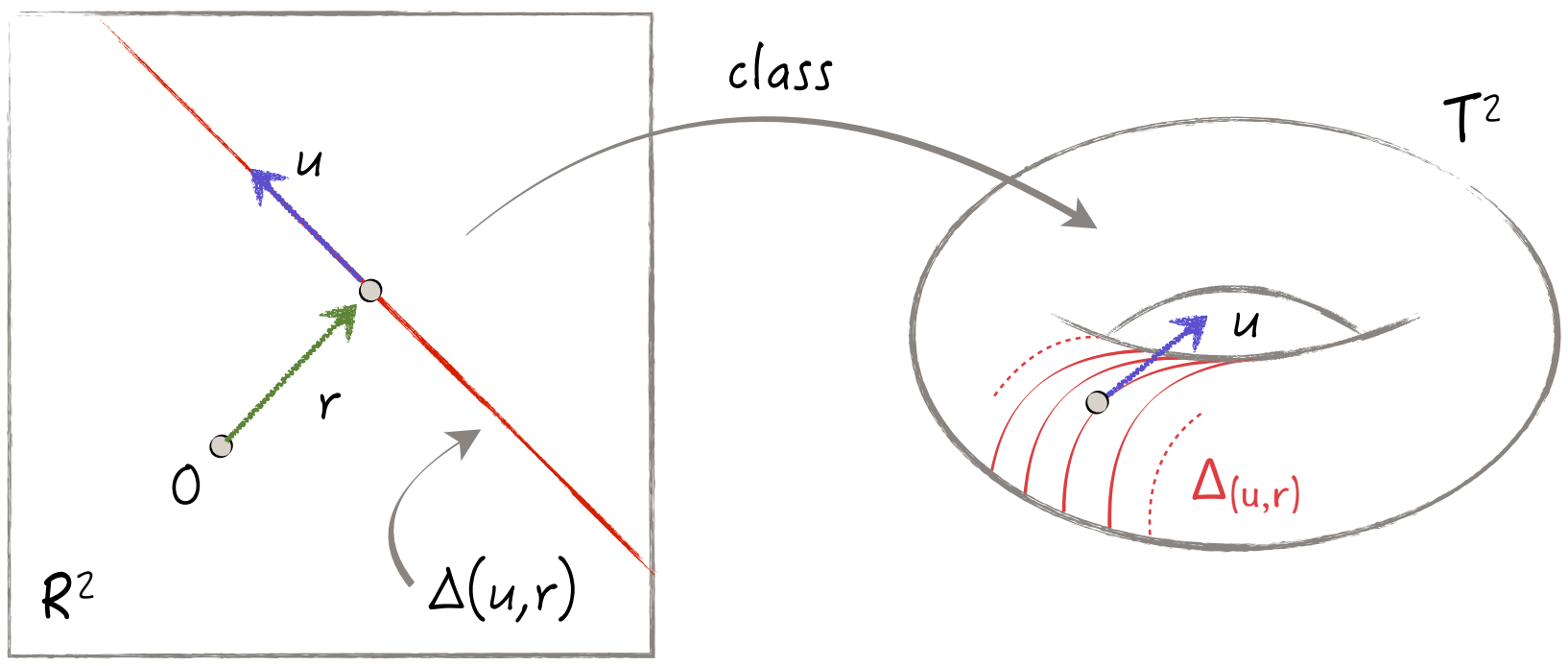

Figure 1. The geodesic trajectories of the torus.

Historically, stratification theory developed a complex axiomatic framework and defined stratifications even for smooth manifolds, despite its original focus on singular subspaces of Euclidean spaces, such as submanifolds with boundaries or corners, or algebraic subspaces. This is no longer essential, as diffeology can directly handle these subspaces, equipped with their subset diffeology, within the framework of diffeological spaces with Klein strata. To illustrate an application of Klein stratification, consider the space of geodesics of the 2-torus, which combines strata that are manifolds with strata that are not.

**Application: The Space of Geodesics of the 2-Torus.** To illustrate an application of Klein stratification, consider the space of geodesics of the 2-torus, which combines strata that are manifolds with strata that are not. This example highlights a curious asymmetry in classical geometry. When we consider the space of oriented geodesic trajectories on the simplest surfaces, we find that for the plane $\mathbf{R}^2$, this space is the cylinder $\mathrm{S}^1 \times \mathbf{R}$, and for the sphere $\mathrm{S}^2$, it is another sphere $\mathrm{S}^2$. Yet for the 2-torus $\mathrm{T}^2$, classical geometry provides no such satisfactory answer; the resulting quotient space is not a manifold and is left without a clear geometric identity. Diffeology repairs this conceptual injustice, endowing the space of geodesics of the torus with its own intrinsic smooth structure and a proper identity within the world of geometric objects.

The geodesics[9] of the torus $\mathrm{T}^2$ are the projections onto the torus of the geodesics of $\mathbf{R}^2$. The set of oriented affine lines in $\mathbf{R}^2$ is diffeomorphic to the tangent space of the circle $\mathrm{S}^1$. An affine line is uniquely determined by a pair $(u,r) \in \mathrm{S}^1 \times \mathbf{R}^2$ such that $u \cdot r = 0$, where the centered dot denotes the inner product. The line associated with the pair $(u,r)$ is $\Delta(u,r) = \{r + tu\}_{t \in \mathbf{R}}$. This space is also equivalent to $\mathrm{S}^1 \times \mathbf{R}$ via $(u,r) \mapsto (u, \rho = u \cdot \mathrm{J}r)$, where J is a rotation by $\pi/2$. The space $\mathrm{Geod}(\mathrm{T}^2)$ of geodesics of $\mathrm{T}^2$ is then diffeomorphic to the quotient $(S^1 \times \mathbf{R})/\mathbf{Z}^2$, where the action of $\mathbf{Z}^2$ is given by $\underline{k}(u,\rho) = (u, \rho + u \cdot k)$ for all $k = (m,n) \in \mathbf{Z}^2$. Letting $u = (\cos(\theta), \sin(\theta))$, the action of $\mathbf{Z}^2$ on $S^1 \times \mathbf{R}$ can be written as

$$\underline{(m,n)}(u,\rho) = (u, \rho + m\cos(\theta) + n\sin(\theta)).$$

[9]We consider oriented geodesic trajectories.

Define

$$\mathrm{class}(u,\rho) \equiv (u, \mathrm{class}_u(\rho)) = (u, \{\rho + k \cdot u \mid k \in \mathbf{Z}^2\}),$$

and let $\mathrm{pr} : (u, \mathrm{class}_u(\rho)) \mapsto u$ be the projection of $\mathrm{Geod}(T^2)$ onto $S^1$. The fibers of this projection are the tori

$$\mathrm{pr} : \mathrm{Geod}(T^2) \to S^1, \text{ with } T_u = \mathrm{pr}^{-1}(u) = \mathbf{R}/[\mathbf{Z}\cos(\theta) + \mathbf{Z}\sin(\theta)].$$

The torus $T_u$ is irrational when $\cos(\theta)$ and $\sin(\theta)$ are linearly independent over $\mathbf{Q}$, and rational, diffeomorphic to a circle of the form $\mathbf{R}/a\mathbf{Z}$, otherwise. Therefore, the space of geodesics of the 2-torus has a diffeology that combines manifolds, i.e., circles when the direction vector $u$ is rational, meaning $\mathbf{R}u$ intersects the lattice $\mathbf{Z}^2$ non-trivially, and pure diffeological irrational tori otherwise, when $u$ is irrational, meaning $\mathbf{R}u$ intersects $\mathbf{Z}^2$ only at the origin.

The Klein stratification of $\mathrm{Geod}(T^2)$ is precisely described in [PIZ25]:

THEOREM. *The orbits of* $\mathrm{Diff}(\mathrm{Geod}(T^2))$ *are the preimages under the projection* pr *of the orbits of* $\mathrm{GL}(2,\mathbf{Z})$ *acting on the circles* $S^1$ *by:* $M(\mathbf{R}u) = \mathbf{R}(Mu)$ *for all* $M \in \mathrm{GL}(2,\mathbf{Z})$.

This example illustrates how the analysis of the internal structure of a relatively complex space, such as the space of geodesics of the 2-torus, can benefit from the diffeological concept of Klein stratification.

**9. Orbifolds and Quasifolds.** As we have seen with manifolds, modeling spaces in diffeology is a powerful tool for defining new subcategories within the category of diffeology [PIZ13, Chap. 4]. Each subcategory starts within diffeology. A manifold is first given with its diffeology, and then its diffeology satisfies some property that makes it a manifold: in this case, to be locally diffeomorphic to a Euclidean space. This approach favors the mathematical intuition that a manifold is something that looks like a Euclidean space, locally, everywhere. We are not burdened, in the definition, by gluing injections according to transition functions. The local diffeomorphisms do the job for us. And many manifolds are indeed built that way, surfaces defined as level sets of smooth functions on Euclidean spaces, like spheres for example, or are quotients of Euclidean spaces, like tori, and the list doesn't stop here.

The concept of orbifolds is of the same kind. First mathematicians encountered spaces that looked like manifolds almost everywhere, except on some subset where they are locally obtained as a quotient by a finite group action. The first example that comes to mind is the quotient $\mathcal{Q}_m = \mathbf{C}/\mathcal{U}_m$ of the complex plane $\mathbf{C}$ by the group $\mathcal{U}_m = \{\exp(2i\pi k/m) \mid k = 0, \dots m-1\}$ of $m$-th roots of unity. A main difficulty starting with topology is that $\mathcal{Q}_m$ is topologically equivalent to $\mathbf{C}$. The singularity at the origin is eaten by the topology. So, even if intuitively an orbifold is something that looks like a quotient of a Euclidean space locally, everywhere, we don't know how to express technically this resemblance, certainly not just involving topology. This is what led Ishiro Satake in his two papers to define what he called *V-manifolds* by a lengthy, complicated, definition [IS56, IS57]. The simple intuition that a V-manifold is just some space locally looking like a quotient by a finite linear group was lost. Later, William Thurston, in his lecture on the subject [Thu78], changed the name of V-manifold to the more sexy *orbifold*, but didn't change the definition.

Here again diffeology gives a simple answer [IKZ10]. Indeed:

DEFINITION (ORBIFOLD). *An orbifold is defined as a diffeological space that is locally diffeomorphic to some quotient* $\mathbf{R}^n/\Gamma$*, where* $\Gamma$ *is a finite subgroup of* $\mathrm{GL}(n,\mathbf{R})$*. The group* $\Gamma$ *may change from point to point.*

Following the same modeling principle, this definition is naturally extended to quasifolds [IZP21], which were originally introduced by Elisa Prato [Pra01]. [10]

DEFINITION (QUASIFOLD). *A quasifold is defined as a diffeological space that is locally diffeomorphic to some quotient* $\mathbf{R}^n/\Gamma$*, where* $\Gamma$ *is a* ***countable*** *subgroup of the affine group* $\mathrm{Aff}(\mathbf{R}^n)$.

This generalization from finite groups to countable affine groups is a significant one. The irrational torus, for example, sits in this subcategory. Almost everything said for orbifolds continues for quasifolds, with the notable exception of the Klein stratification, which is no longer necessarily standard (i.e., its strata may not be manifolds). We have already seen an example with $\mathrm{Geod}(\mathrm{T}^2)$, a quasifold which contains strata that are themselves irrational tori.

But that definition would be empty if there was not the following theorem (op. cit.), where intuition reclaims its rights:

THEOREM. *Every Satake V-manifold is a diffeological orbifold for a natural diffeology associated with a Satake's* defining family. *Conversely, to every diffeological orbifold*[11] *is associated a Satake's defining family that is a V-manifold. These operations are inverse to each other.*

To illustrate how diffeology handles this definition, consider the description of a plot $\zeta$ of the teardrop orbifold structure on the 2-sphere $\mathrm{S}^2 \subset \mathbf{C} \times \mathbf{R}$. Let N denote the north pole:

(A) If $\zeta(r_0) \neq \mathrm{N}$, then there exists a ball B centered at $r_0$ such that $\zeta \restriction \mathrm{B}$ is smooth.

(B) If $\zeta(r_0) = \mathrm{N}$, there exist a ball B centered at $r_0$ and a smooth parametrization $z$ in $\mathbf{C}$ defined on B such that, for all $r \in \mathrm{B}$,

$$\zeta(r) = \frac{1}{\sqrt{1+|z(r)|^{2m}}}\begin{pmatrix} z(r)^m \\ 1 \end{pmatrix}.$$

Here, the integer $m > 1$ determines the order of the singularity at the north pole where the space is locally diffeomorphic to $\mathbf{C}/\mathcal{U}_m$.

---

[10] Note that the diffeological definition of a quasifold, based on local models, differs from Prato's original construction. This adaptation follows the same philosophy used for orbifolds: we take the intuitive local picture ('quotient by a group action') as the defining characteristic of an existing diffeology, in order to integrate the concept seamlessly into the modeling framework. This parallels how the diffeological definition of an orbifold differs from Satake's original axiomatic approach.

[11] Modulo a restriction imposed in Satake definition, on the absence of reflexions.

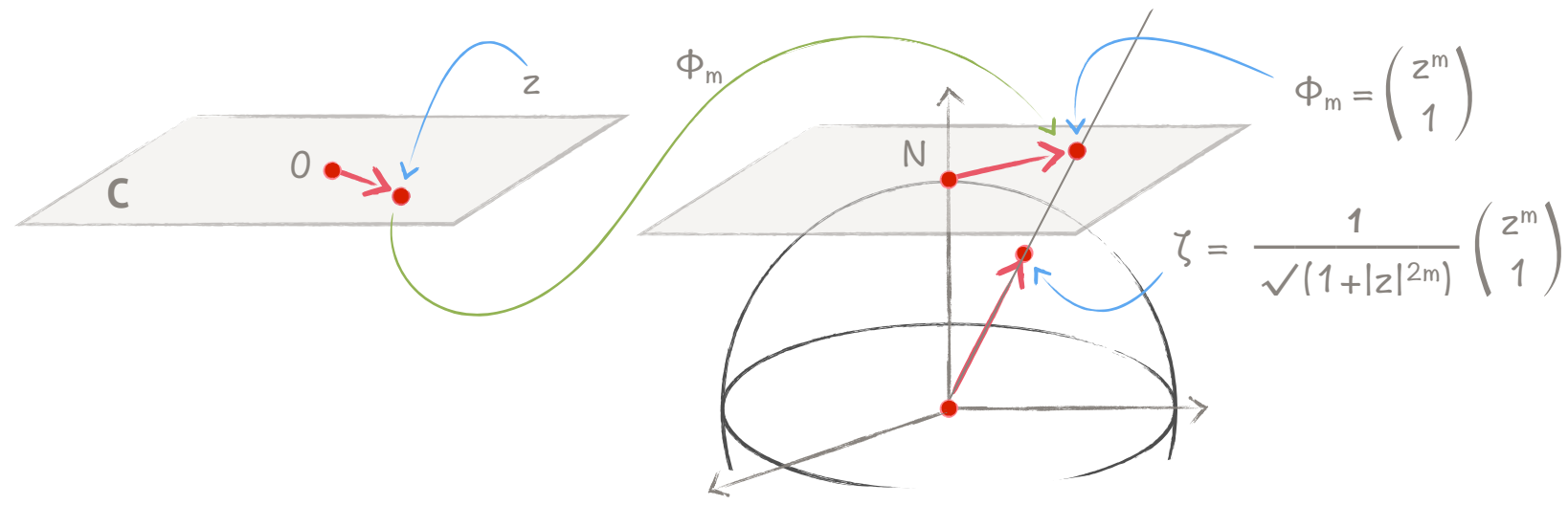

Figure 2. The TearDrop Orbifold structure on the Sphere $S^2$.

By embedding orbifolds into the category {Diffeology} we also gain the morphisms between orbifolds, which are naturally smooth maps in the sense of the diffeology. This may seem like a small gain, but it was a real difficulty, since Satake was unable to give a satisfactory notion of smooth maps between orbifolds. This is what he wrote himself in [IS57, page 469] as a footnote:

> "*The notion of* $C^\infty$ *-map thus defined is inconvenient in the point that a composite of two* $C^\infty$ *-maps defined in a different choice of defining families is not always a* $C^\infty$ *map.*"

Indeed, it was too much to ask for a smooth map to lift locally on the local model equivariantly with respect to the structure group $\Gamma$. There is a famous counter-example elaborated in [IKZ10, Ex. 24] of a smooth map from simple orbifolds that has no locally equivariant lifting.

Some people may argue that this construction is too simplistic and evacuate deeper concepts related to another definition, more complicated, by way of groupoids. This would be an argument if it was impossible to find, starting with the geometry of the orbifold, that is, its diffeology, a groupoid capturing this complexity, but this has been done in [IZL18] by associating categorically a *structure groupoid* to every diffeological orbifold.

*The Structure Groupoid.* The construction begins by choosing a presentation of the orbifold X, that is,

(1) The local model is defined by *charts* which are local diffeomorphism $f$ from $\mathbf{R}^n/\Gamma$ to X.

(2) Then, by gathering the charts in an *atlas*, that is a family $\mathscr{A}$ of charts $f : U \to X$ whose values $f(U)$ cover X.

(3) Associated to this atlas is the *strict generating family* $\mathscr{F}$ made by the plots $F = f \circ \mathrm{class}$, defined on $\tilde{U} = \mathrm{class}^{-1}(U)$, with class : $\mathbf{R}^n \to \mathbf{R}^n/\Gamma$.

(4) The *nebula* $\mathscr{N}$ of $\mathscr{F}$ is then defined as the sum of the domains of F, and the evaluation map ev : $\mathscr{N} \to X$.

$$\mathscr{N} = \coprod_{F \in \mathscr{F}} \mathrm{dom}(F) \quad \text{and} \quad \mathrm{ev} : (F, r) \mapsto F(r)$$

Actually, X identifies to the quotient of $\mathscr{N}$ by ev. In technical terms, ev is a *subduction.*

(5) Now, define the groupoid $\mathbf{G}$ whose set of objects is the nebula $\mathcal{N}$, whose morphisms between $(\mathrm{F}, r)$ and $(\mathrm{F}', r')$ is not empty only if $\mathrm{F}(r) = \mathrm{F}'(r')$, and are the germs of the local diffeomorphisms $\varphi$ of $\mathcal{N}$ that project onto the identity of X by the evaluation map ev, that is, such that $\mathrm{F} =_{\mathrm{loc}} \mathrm{F}' \circ \varphi$.

DEFINITION (STRUCTURE GROUPOID). *This groupoid* $\mathbf{G}$ *is the* structure groupoid *associated to the atlas* $\mathcal{A}$. *It is equipped with a special diffeology defined as the quotient by germification of some functional diffeology on local smooth maps, and becomes a diffeological groupoid.*

THEOREM (EQUIVALENCE OF STRUCTURES GROUPOIDS). *The orbifold* X *identifies with the space of transitivity components*[12] *of* $\mathbf{G}$,

$$\mathrm{X} \equiv \mathcal{N}/\mathrm{Transitivity}.$$

*Different atlases give equivalent structure groupoids. Diffeomorphic quasifolds have equivalent structure groupoids.*

*Sketch of Proof.* This the consequence of a fundamental theorem proved in [IKZ10] that every local diffeomorphism of a groupoid lifts locally on $\mathcal{N}$ as a local equivariant diffeomorphism. And it is extended to quasifold in [IZP21]. It is important to note that this is not the case for arbitrary maps, as show the example [IKZ10, Ex. 24]. ►

Thus, the equivalence class[13] of a structure groupoid is a diffeological invariant that is naturally associated with the orbifold. Any of these groupoids $\mathbf{G}$ contains exactly the same information as the *effective orbifold groupoids* found in mathematical literature, see [KM23].

The non-effective part is a structure over the groupoid, which may add something to the groupoid but doesn't change its space of transitivity components — the quotient space of by the transitivity relationship — that we call the base *base* space of the groupoid, viewing the groupoid itself as a structure built over the underlying quasifold. This is analogous to the situation with fiber bundles over a fixed **base space**: many different fiber bundles can be constructed over the same base, but these additional structures do not alter the intrinsic geometry of the base space itself. However, the category of diffeological groupoids $\mathbf{G}$ over the same base $\mathrm{X} = \mathrm{Obj}(\mathbf{G})/\mathrm{Transitivity}$ could be an interesting subject of study.

The structure groupoid that characterizes the internal smooth structure of the orbifold reflects on its Klein stratification. Thanks to the following fact, proved in [GIZ23]:

PROPOSITION. *Two points $x$ and $x'$ in an orbifold* X *are on the same Klein stratum if and only if their isotropies are conjugate modulo* $\mathrm{GL}(n, \mathbf{R})$.

In other words, the Klein stratification of the orbifold, which expresses the internal singularities of the space, follows the structure groupoid. Moreover, it is proved in this paper that:

---

[12] In much of the standard literature on groupoids, these equivalence classes are referred to as *orbits*. We have chosen the term *transitivity components* in this context to avoid potential ambiguity. Specifically, we wish to reserve the term 'orbit' for the action of a group or pseudogroup, such as the orbits of local diffeomorphisms that constitute the Klein strata of the resulting quotient space.

[13] Here, 'equivalence' is taken in the standard sense of category theory [McL71, Chap. IV]. We will see later another kind of equivalence leading to the Morita equivalence of associated C*-algebras.

THEOREM (GIZ). *The Klein stratification of a (locally finite) diffeological orbifold is a standard stratification: the strata are locally closed manifolds, and the set of strata, equipped with the quotient diffeology, satisfies the* $T_0$ *separation condition.*

**10. Application: The Definitive Geometry of Orbit Spaces.** The study of orbit spaces $M/G$ is a cornerstone of geometry, yet a classical problem has persisted: the standard tool for analyzing the action, the *orbit-type stratification* on $M$, is an extrinsic structure that fails to reliably predict the *intrinsic* geometry of the singular quotient $M/G$. Our recent work [GIZ25] provides the definitive solution to this problem within the diffeological framework.

The intrinsic geometry of the quotient is captured by its **Klein stratification**, the partition into orbits of local diffeomorphisms. The correct predictive tool on the manifold M is not the orbit-type stratification, but the finer **isostabilizer decomposition**, whose elements are the connected components of submanifolds where the stabilizer subgroup is constant. The central result is the **Correspondence Theorem**, which establishes a surjective map from the isostabilizer decomposition of M to the Klein stratification of $M/G$. This provides the exact dictionary between the action's local data and the quotient's singularity types.

This framework yields a concrete, computable invariant — the **diffeological dimension map** — which satisfies the formula:

$$\dim_{\mathscr{O}}(M/G) = \dim(M) - \dim(\mathscr{O}).$$

This formula quantitatively captures how the "collapsing" of orbits of different dimensions creates singularities of varying diffeological dimensions in the quotient.

Most profoundly, this analysis reveals that the classical orbit-type stratification is fundamentally **incomparable** to the intrinsic Klein stratification of the quotient. The relationship is not a simple matter of refinement, as two key examples demonstrate [PIZ91, Aud91].

First, in the action of SO(3) on the complex projective plane $P^2(\mathbf{C})$, two distinct singular orbit-type strata — one with stabilizer type (SO(2)) and another with type (O(2)) — collapse into a *single* Klein stratum in the quotient. Here, the classical stratification is too fine.

Conversely, in a specific action of $U(1)$ on $P^2(\mathbf{C})$, a *single* orbit-type stratum — the set of three fixed points — is torn apart by the projection. Two of the fixed points map to one Klein stratum, while the third, having a non-equivalent slice representation, maps to a completely different one. Here, the classical stratification is too coarse.

This leads to a fundamental conclusion: the quotient space $M/G$ possesses an intrinsic smooth identity, which we term an **orthofold**, that can be independent of the specific global action that created it. The quotient "forgets" the details of the orbit-type structure, remembering only the local diffeomorphic types of its singularities. This work thus replaces a classical, often inadequate tool with a precise and predictive theory for one of geometry's most fundamental constructions.

## 3. Symplectic Geometry in the Diffeological Setting

As is well-known, the treatment of singularities is a recurring challenge in differential geometry and exterior differential calculus. Two principal areas where these challenges arise are: singular quotients and infinite-dimensional spaces, both of which are particularly pertinent in symplectic geometry.

Singularities are especially critical in the context of symplectic reduction. Symplectic reduction involves reducing a coisotropic submanifold W within a symplectic manifold M by the kernel of the restricted symplectic form $\omega$. This process was thoroughly addressed in the seminal work of Alan Weinstein and Jerrold Marsden [MW74] for the case where the quotient space $Q = W/\ker(\omega \upharpoonright W)$ is a manifold. However, the case where Q is merely a topological space, possibly non-Hausdorff, remained largely unaddressed. In recent decades, singular reduction has been addressed through staged approaches along stratifications, analogous to Marsden-Weinstein reduction, as seen in works like [LS91], and [Wat12] (which employs diffeological techniques).

Diffeology offers an original solution to this problem, which should be addressed in its full generality. Since exterior differential calculus extends naturally to diffeological spaces, and quotient spaces are always equipped with the quotient diffeology, the question can be reformulated as

> *Question:* Let $Q = W/\ker(\omega \upharpoonright W)$ be equipped with the quotient diffeology. Does there exist a 2-form $\tilde{\omega}$ on Q such that $\omega \upharpoonright W = \mathrm{class}^*(\tilde{\omega})$?

Because manifolds form a full and faithful subcategory of {Diffeology}, when Q is a manifold, then $\tilde{\omega}$ coincides with the Marsden-Weinstein construction. A complete answer to this question is not yet available, but we have partial results. We will examine two examples in the case of a 1-dimensional kernel.

The second problem that extends symplectic geometry beyond its traditional boundaries emerged with the development of classical and quantum field theories. These theories introduced symplectic-type structures on infinite-dimensional spaces such as spaces of functions and connections. We will explore how tools like the moment map generalize to diffeology and apply to these cases. In this section, we will also consider an example of symplectic reduction in an infinite-dimensional setting with singularities. The combination of these two challenging situations further demonstrates the versatility of diffeology.

**11. Differential Forms and Exterior Calculus.** A *differential $k$-form* $\alpha$ on a diffeological space X assigns to each plot P in X a smooth form

$$\alpha(\mathrm{P}) \in \Omega^k(\mathrm{dom}(\mathrm{P})),$$

such that, for all smooth parametrizations F in dom(P),

$$\alpha(\mathrm{P} \circ \mathrm{F}) = \mathrm{F}^*\big(\alpha(\mathrm{P})\big).$$

The set $\Omega^k(\mathrm{X})$ of differential $k$-forms on X is a real vector space.

Let $f : \mathrm{X} \to \mathrm{X}'$ be a smooth map, and $\alpha' \in \Omega^k(\mathrm{X}')$. The *pullback* of $\alpha'$ by $f$ is defined by:

$$f^*(\alpha')(\mathrm{P}) = \alpha'(f \circ \mathrm{P}).$$

Note that the value of a form $\alpha$ on a plot P is simply the pullback of $\alpha$ by P: $\alpha(\mathrm{P}) = \mathrm{P}^*(\alpha)$.

Exterior differentiation is defined by the operator:

$$d : \Omega^k(\mathrm{X}) \to \Omega^{k+1}(\mathrm{X}), \ \text{ given by } \ d\alpha \colon \mathrm{P} \mapsto d[\alpha(\mathrm{P})]$$

The operator $d$ is linear and smooth with respect to the functional diffeology. It satisfies $d \circ d = 0$.

Crucially, the following criterion determines when a differential form descends to a quotient. Consider a subduction $\pi : \mathrm{Y} \to \mathrm{X}$ between diffeological spaces, meaning that X is equivalent to the quotient $\mathrm{Y}/\pi$.

CRITERION. *Let $\alpha \in \Omega^k(\mathrm{Y})$. There exists a differential form $\beta \in \Omega^k(\mathrm{X})$ such that $\alpha = \pi^*(\beta)$ if and only if: for any two plots $\mathrm{P}, \mathrm{P}'$ in $\mathrm{Y}$, if $\pi \circ \mathrm{P} = \pi \circ \mathrm{P}'$, then $\alpha(\mathrm{P}) = \alpha(\mathrm{P}')$.*

This criterion is a fundamental tool for determining if a differential form descends to the quotient by an equivalence relation. This is essentially all we need for the subsequent discussion.

**12. Reduction by the Geodesic Flow.** The space of geodesics[14] of a Riemannian manifold $(\mathrm{M}, g)$ is a particular case of symplectic reduction that may be singular.

Consider the 1-form on the tangent space TM

$$\lambda = \bar{u}\, dx, \ \text{ or in a local chart: } \ \lambda = \sum_{i=1}^{n} u_i\, dx^i, \ \text{ with } \ u_i = g_{ij} u^j.$$

It is the pullback of the canonical Liouville 1-form on T*M by the metric transposition $u \mapsto \bar{u}$, with $\bar{u}(v) = g(u, v)$. Its differential $d\lambda$ is a symplectic form, and its restriction to the unit tangent bundle UM, the subspace of tangent vectors of norm 1, is a contact form:

$$\dim(\ker d\lambda) = 1 \ \text{ and } \ \ker(\lambda) \cap \ker(d\lambda) = \{0\}.$$

The kernel of the restriction of $d\lambda$ on UM is generated by the Reeb vector $\xi$: $\ker(d\lambda \restriction \mathrm{UM}) = \mathbf{R}\xi$ and $\lambda(\xi) = 1$. The integral curves of this linear distribution project onto M as geodesic trajectories. The set of characteristics of $d\lambda \restriction \mathrm{UM}$ thus identifies with the set of geodesics of $g$ on M:

$$\mathrm{Geod}(\mathrm{M}) \equiv \mathrm{UM}/\ker(d\lambda).$$

When Geod(M) is a manifold for the quotient diffeology, $d\lambda \restriction \mathrm{UM}$ descends to Geod(M) as a symplectic form $\omega$. This is the case, for example, for the sphere $\mathrm{S}^2$, where $\mathrm{Geod}(\mathrm{S}^2) \equiv \mathrm{S}^2$ via the mapping $(x, u) \mapsto \ell = x \wedge u$.

The behavior of Geod(M) when it is not a manifold is beyond the scope of ordinary differential geometry, but diffeology elucidates it as a consequence of this theorem:

[14] As before, the term "geodesic" refers to an oriented *geodesic trajectory*, which we distinguish from a parameterized *geodesic curve*. This distinction, while always present, has profound structural implications that become especially critical in the pseudo-Riemannian setting. Indeed, focusing on the space of trajectories provides the key to a unified geometric description that classical reduction methods fail to achieve, a topic explored in detail in [PIZ25c].

THEOREM (REDUCTION OF A CONTACT FORM). *Let* $\mathrm{Y}$ *be a Hausdorff and second-countable manifold equipped with a contact form* $\lambda$*. There always exists on the space* $\mathscr{C}$ *of characteristics of* $d\lambda$ *a parasymplectic form* $\omega$*, symplectically generated, such that*

$$\mathrm{class}^*(\omega) = d\lambda,$$

*where* $\mathrm{class} : \mathrm{Y} \to \mathscr{C} = \mathrm{Y}/\ker(d\lambda)$ *is the projection. If* $\mathrm{Y}$ *has dimension* $2n-1$ *(contact manifolds are odd dimensional), then* $\mathscr{C}$ *has dimension* $2n-2$*. When* $\mathscr{C}$ *is a manifold, then* $\omega$ *is symplectic.*

Because it is important to understand why such a theorem holds in diffeology, I will provide a sketch of the proof that explains how it works. For more details, see [PIZ25, p.275]

*Sketch of Proof.* The proof relies on the properties of the Reeb vector field $\xi$, whose flow preserves $\lambda$: $\exp(t\xi)^*(\lambda) = \lambda$. It consists in proving that $d\lambda$ satisfies the condition to descend to the quotient space $\mathscr{C}$. Let $\mathrm{P}$ and $\mathrm{P}'$ be two plots in $\mathrm{Y}$ such that $\mathrm{class} \circ \mathrm{P} = \mathrm{class} \circ \mathrm{P}'$. This means that for any parameter value $r$, $\mathrm{P}(r)$ and $\mathrm{P}'(r)$ lie on the same characteristic curve of $d\lambda$. This implies the existence of a smooth parametrization $s(r)$ such that $\mathrm{P}'(r) = \exp(s(r)\xi)(\mathrm{P}(r))$, thanks to a standard result in the theory of ordinary differential equations. Using the fact that the flow preserves $\lambda$, we have: $\lambda(\mathrm{P}')_r = \lambda(\exp(s(r)\xi)(\mathrm{P}(r)))_r = \lambda(\mathrm{P})_r + d s_r$. Taking the exterior derivative of both sides, we get $d\lambda(\mathrm{P}') = d\lambda(\mathrm{P})$. Therefore, by the criterion for differential forms to descend to the quotient in diffeology, there exists a (closed) 2-form $\omega$ on $\mathscr{C}$ such that $\mathrm{class}^*(\omega) = d\lambda$. ►

COROLLARY. *As a reduction of the unit tangent bundle* $\mathrm{UM}$ *by the kernel of* $\lambda$*, the space of geodesics of a Riemannian manifold is always a parasymplectic space for the projection* $\omega$ *of* $d\lambda \restriction \mathrm{UM}$.

This applies in particular to the geodesics of the 2-torus, for which, as we have seen, $\mathrm{Geod}(\mathrm{T}^2)$ is not a manifold but combines circles and irrational tori depending on the slopes. Indeed, $\mathrm{Geod}(\mathrm{T}^2)$ is the quotient space by the action of the additive group $\mathbf{R}$ whose stabilizers vary from point to point, ranging from $\{0\}$ to the groups $\mathscr{U}_m$ of the $m$th roots of unity, where $m$ runs through $\mathbf{N}$. This clearly does not satisfy the conditions of the Marsden-Weinstein procedure, and the quotient space is even non-Hausdorff. The singularities of this space, which are the different orbits of the local diffeomorphisms, do not prevent the parasymplectic form $\omega$ from being well-defined on $\mathrm{Geod}(\mathrm{T}^2)$.

This is an example of symplectic reduction in the presence of singularities.[15] It demonstrates how diffeology allows us to extend the concept of symplectic construction to spaces that are not manifolds. The space of geodesics on the torus, despite its non-manifold nature, carries a natural parasymplectic structure inherited from the contact form on its unit tangent bundle. The ability to work with such structures on these types of spaces highlights the power and versatility of the diffeological framework and unifies a geometric construction: the parasymplectic character of the space of geodesics, which is otherwise impossible.

[15]The precise definition of a 'symplectic' diffeological space is a subject of ongoing discussion. One approach, which we favor, defines a space as symplectic if it satisfies a generalized Darboux condition: the transitivity of its pseudogroup of local symplectomorphisms and the injectivity of its universal moment map. This definition, however, excludes some objects like standard symplectic orbifolds, where the local symmetry group is not transitive due to the singularities. An alternative approach defines a space as symplectic if its diffeology can be generated by plots on which the pullback of the 2-form is symplectic. We refer to spaces satisfying this latter condition as 'symplectically generated'.

**13. Reduction of Infinite Dimensional Ellipsoids.** We shall now examine an example of symplectic reduction in diffeology within an infinite-dimensional space where the reduction is singular. This example combines two situations that standard symplectic geometry cannot readily handle. It demonstrates that in this infinite-dimensional case, despite the singularities, the symplectic form descends to the reduced space when restricted to a generic level of the Hamiltonian. For detailed technical information, refer to [PIZ16].

Consider the space of periodic functions from $\mathbf{R}$ to $\mathbf{C}$:

$$\mathscr{C}^\infty_{\mathrm{per}}(\mathbf{R},\mathbf{C}) = \{f \in \mathscr{C}^\infty(\mathbf{R},\mathbf{C}) \mid f(x+1) = f(x)\}.$$

Let Surf be the standard symplectic form on $\mathbf{C}$,

$$\mathrm{Surf}_z(\delta z, \delta' z) = \frac{1}{2i}\left[\delta\bar z\, \delta' z - \delta'\bar z\, \delta z\right],$$

where $z, \delta z, \delta' z \in \mathbf{C}$. For all $x \in \mathbf{R}$, let

$$\hat x : \mathscr{C}^\infty_{\mathrm{per}}(\mathbf{R},\mathbf{C}) \to \mathbf{C} \quad \text{with} \quad \hat x(f) = f(x)$$

be the *evaluation map*. The pullback $\hat x^*(\mathrm{Surf})/\pi$ is a closed 2-form on $\mathscr{C}^\infty_{\mathrm{per}}(\mathbf{R},\mathbf{C})$. Let $\omega$ be its mean value:

$$\omega = \frac{1}{\pi}\int_0^1 \hat x^*(\mathrm{Surf})\, dx.$$

Evaluated on a plot $\mathrm{P} : r \mapsto f_r$, we obtain:

$$\omega(\mathrm{P})_r(\delta r, \delta' r) = \frac{1}{2i\pi}\int_0^1 \left\{ \frac{\partial \overline{f_r(x)}}{\partial r}(\delta r)\frac{\partial f_r(x)}{\partial r}(\delta' r) - \frac{\partial \overline{f_r(x)}}{\partial r}(\delta' r)\frac{\partial f_r(x)}{\partial r}(\delta r)\right\} dx.$$

The vector space $\mathscr{C}^\infty_{\mathrm{per}}(\mathbf{R},\mathbf{C})$ equipped with $\omega$ is a *diffeological symplectic space* in a strong sense.[16]

We then transport this structure to the space $\mathscr{E}$ of Fourier coefficients via $j : f \mapsto (f_n)_{n\in\mathbf{Z}}$. A parametrization $\mathrm{P} : r \mapsto (f_n(r))_{n\in\mathbf{Z}}$ in $\mathscr{E}$ is a plot if:

(1) The functions $f_n : \mathrm{dom}(\mathrm{P}) \to \mathbf{C}$ are smooth.
(2) Their successive derivatives are uniformly rapidly decreasing:

$$n^p \frac{\partial^k f_n(r)}{\partial r^k} \xrightarrow[|n|\to\infty]{} 0.$$

Next, the infinite product of tori $\mathrm{T}^\infty = \prod_{n\in\mathbf{Z}} \mathrm{T}$ is equipped with a *tempered diffeology*. A parametrization $\zeta : r \mapsto (z_n(r))_{n\in\mathbf{Z}}$ in $\mathrm{T}^\infty$ is *tempered* if the $z_n$ are smooth and if for every $k \in \mathbf{N}$, for every $r_0$ in the domain of the parametrization, there exist a closed ball $\bar{\mathrm{B}} \subset \mathrm{dom}(\zeta)$ centered at $r_0$, a polynomial $\mathrm{P}_k$ and an integer N such that

$$\text{For all } r \in \mathrm{B}, \text{ for all } n > \mathrm{N}, \quad \left| \frac{\partial^k z_n(r)}{\partial r^k} \right| \leq \mathrm{P}_k(n).$$

[16]It is an (additive) diffeological group and an (affine) coadjoint orbit of itself.

Equipped with tempered parametrizations, $\mathrm{T}^\infty$ is a diffeological group acting smoothly, by automorphisms, on the symplectic space $\mathscr{E}$ of Fourier coefficients by pointwise multiplication:

$$\text{For all } \tau = (\tau_n)_{n\in\mathbf{Z}} \in \mathrm{T}^\infty \text{ and } \mathrm{Z} = (\mathrm{Z}_n)_{n\in\mathbf{Z}} \in \mathscr{E} : \ \tau\mathrm{Z} = (\tau_n \mathrm{Z}_n)_{n\in\mathbf{Z}}.$$

The moment map, which will be defined later in full generality, of this action is

$$\mu(\mathrm{Z}) = \frac{1}{2i\pi} \sum_{n\in\mathbf{Z}} |\mathrm{Z}_n|^2 \, \pi_n^*(\theta) + \sigma,$$

where $\pi_n$ is the projection on the $n$-th component, $\theta \in \Omega^1(\mathrm{T})$ is the canonical volume form, and $\sigma$ is a constant.

We choose a sequence of rationally independent real numbers $\alpha = (\alpha_n)_{n\in\mathbf{Z}}$, with all $0 < \alpha_n < 1$. We induce the real line in the infinite torus as an irrational solenoid $\mathscr{S}_\alpha = \{(\exp(2i\pi\alpha_n t))_{n\in\mathbf{Z}} \mid t \in \mathbf{R}\}$. The restriction to this additive action of $\mathbf{R}$ has for moment map $\nu$ defined by the Hamiltonian $h$:

$$\nu(\mathrm{Z}) = h(\mathrm{Z})\, dt \ \text{ with } \ h(\mathrm{Z}) = \sum_{n\in\mathbf{Z}} \alpha_n |\mathrm{Z}_n|^2 + c,$$

with $c$ a real constant. Let

$$\mathrm{S}^\infty_\alpha = \left\{ \mathrm{Z} = (\mathrm{Z}_n)_{n\in\mathbf{Z}} \in \mathscr{E} \;\middle|\; \sum_{n\in\mathbf{Z}} \alpha_n |\mathrm{Z}_n|^2 = 1 \right\}.$$

be a generic level subspace of the Hamiltonian of the solenoid action and then reduce it by the solenoid action. Note that $\mathrm{S}^\infty_\alpha$ is an infinite dimensional ellipsoid with axes $\alpha_n$. This action is not free and generates infinitely many singular orbits:

$$\mathrm{S}^1_n = \{\mathrm{Z} \in \mathrm{S}^\infty_\alpha \mid \mathrm{Z}_m = 0 \text{ if } m \neq n\}, \text{ with } n \in \mathbf{Z},$$

isomorphic to circles of radius $1/\sqrt{\alpha_n}$. The reduced space is an *infinite quasi-projective space* (an infinite dimensional quasifold, although this has not yet been rigorously defined):

$$\mathrm{P}^\infty_\alpha = \mathrm{S}^\infty_\alpha / \mathbf{R}.$$

Despite the presence of singular orbits,

THEOREM. *The restriction of the symplectic form $\omega$ on the infinite ellipsoid $\mathrm{S}^\infty_\alpha$ descends to the quotient $\mathrm{P}^\infty_\alpha$ as a closed 2-form $\varpi$, equipping the infinite quasi-projective space with a parasymplectic structure.*

Here, diffeology provides a direct approach to handling infinite dimensions with singularities, avoiding the complexities of functional analysis, which may not be suitable for this problem. The proof, which is entirely diffeological, is outlined below, as we cannot rely on a classical theorem like the one used in the reduction of a contact form.

*Sketch of Proof.* We apply the general criterion mentioned earlier. Let $\mathrm{P} : \mathrm{U} \to \mathrm{S}^\infty_\alpha$ and $\mathrm{P}' : \mathrm{U} \to \mathrm{S}^\infty_\alpha$ be two plots such that $\mathrm{class} \circ \mathrm{P} = \mathrm{class} \circ \mathrm{P}'$. We consider separately the subset $\mathscr{S} = \cup_{n\in\mathbf{Z}} \mathrm{S}^1_n$ of singular orbits and the complementary set of principal orbits $\mathrm{S}^\infty_\alpha - \mathscr{S}$ (where the action is free). It is a fact that $\mathrm{S}^\infty_\alpha - \mathscr{S}$ is open in the D-topology. Thus, $\mathrm{U}_0 = \mathrm{P}^{-1}(\mathrm{S}^\infty_\alpha - \mathscr{S})$ is an open subset of $\mathrm{U}$. Since $\mathrm{class} \circ \mathrm{P} = \mathrm{class} \circ \mathrm{P}'$, $\mathrm{P}^{-1}(\mathrm{S}^\infty_\alpha - \mathscr{S}) = \mathrm{P}'^{-1}(\mathrm{S}^\infty_\alpha - \mathscr{S}) = \mathrm{U}_0$.

Now, the restrictions of $\mathrm{P}$ and $\mathrm{P}'$ to $\mathrm{U}_0$ take values in the subset of $\mathrm{S}^\infty_\alpha$ consisting of principal orbits of $\mathbf{R}$, for which the action of $\mathbf{R}$ is free. Let us denote the component functions of the plots by $\mathrm{P}(r) = (\mathrm{Z}_n(r))_{n\in\mathbf{Z}}$ and $\mathrm{P}'(r) = (\mathrm{Z}'_n(r))_{n\in\mathbf{Z}}$. Therefore, for each $r \in \mathrm{U}_0$, there is a unique $\tau(r) \in \mathbf{R}$ such that, for all $n$, $\mathrm{Z}'_n(r) = e^{2i\pi\alpha_n\tau(r)}\mathrm{Z}_n(r)$. Since for all $r_0 \in \mathrm{U}_0$, there exists $n \in \mathbf{Z}$ such that $\mathrm{Z}_n(r_0) \neq 0$, there exists a neighborhood of $r_0$ where $\mathrm{Z}_n(r) \neq 0$. On this neighborhood, $\mathrm{Z}'_n(r) \neq 0$, and $e^{2i\pi\alpha_n\tau(r)} = \mathrm{Z}'_n(r)/\mathrm{Z}_n(r)$. Because $r \mapsto \mathrm{Z}'_n(r)$ and $r \mapsto \mathrm{Z}_n(r)$ are smooth, it follows that $r \mapsto e^{2i\pi\alpha_n\tau(r)}$, and hence $\tau$, are smooth.

Now, $\omega = d\varepsilon$ with:

$$\begin{aligned}\varepsilon(\mathrm{P}')_r(\delta r) = \frac{1}{2i\pi}\sum_{n\in\mathbf{Z}} \overline{\mathrm{Z}'_n(r)}\frac{\partial \mathrm{Z}'_n(r)}{\partial r}(\delta r) &= \frac{1}{2i\pi}\sum_{n\in\mathbf{Z}} \overline{\mathrm{Z}_n(r)}\frac{\partial \mathrm{Z}_n(r)}{\partial r}(\delta r) \\ &\quad + \left(\sum_{n\in\mathbf{Z}} \alpha_n |\mathrm{Z}_n(r)|^2\right)\frac{\partial \tau(r)}{\partial r}(\delta r) = \varepsilon(\mathrm{P})_r(\delta r) + \tau^*(dt)_r(\delta r).\end{aligned}$$

Therefore, $[\omega(\mathrm{P}') - \omega(\mathrm{P})] \restriction \mathrm{U}_0 = 0$. By continuity, $[\omega(\mathrm{P}') - \omega(\mathrm{P})] \restriction \overline{\mathrm{U}}_0 = 0$, where $\overline{\mathrm{U}}_0$ is the closure of $\mathrm{U}_0$.

It remains to examine the complementary subset $\mathrm{V} = \mathrm{U} - \overline{\mathrm{U}}_0$. The subset $\mathrm{V}$ is open, so $\mathrm{P} \restriction \mathrm{V}$ and $\mathrm{P}' \restriction \mathrm{V}$ are plots of $\mathrm{S}^\infty_\alpha$ with values in the subset of singular orbits $\mathscr{S}$. Since $\mathscr{S}$ is a union of disjoint 1-dimensional submanifolds (the circles $\mathrm{S}^1_n$), the image of any plot with domain in $\mathrm{V}$ must lie within this 1-dimensional set. Consequently, the evaluation of the 2-form $\omega$ on any such plot must be zero. Therefore, the restrictions $\omega(\mathrm{P} \restriction \mathrm{V})$ and $\omega(\mathrm{P}' \restriction \mathrm{V})$ are both zero. In conclusion, $\omega(\mathrm{P}') = \omega(\mathrm{P})$ everywhere on $\mathrm{U}$. This proves that there exists a 2-form $\varpi$ on $\mathrm{P}^\infty_\alpha = \mathrm{S}^\infty_\alpha/\mathbf{R}$ such that $\mathrm{class}^*(\varpi) = \omega$. ►

The restriction of the symplectic form $\omega$ on the infinite ellipsoid $\mathrm{S}^\infty_\alpha$ descends to the quotient $\mathrm{P}^\infty_\alpha$ as a closed 2-form $\varpi$, equipping the infinite quasi-projective space with a parasymplectic structure. This establishes that the reduction process yields a meaningful geometric structure in a challenging infinite-dimensional and singular setting, providing a foundation for further investigation of these geometric objects within diffeology. We anticipate analogous results for general symplectic actions of Lie groups. The case of symplectic reduction for general symplectic diffeological spaces may present different challenges.

## 4. Cartan-de Rham Calculus and Homotopy Invariance

Not only does diffeology provide solutions to conceptual problems such as the treatment of singularities, but it also offers elegant alternatives to sometimes tedious theorem proofs, as exemplified by the homotopy invariance of de Rham cohomology. This invariance states that the pullback of a closed differential form by homotopic maps yields cohomologous forms. Consequently, the morphism between de Rham cohomology groups induced by the pullback of smooth maps depends solely on the homotopy class of the maps.

To prove this invariance within diffeology, we introduce the *Chain Homotopy Operator* ($\mathsf{K}$), a crucial tool that offers a direct and elegant proof of this theorem for diffeological spaces, and hence for manifolds, given that they form a full and faithful subcategory of {Diffeology}. Diffeology's contribution in this context, as we shall see, lies in enabling us to work within the category with the space of paths and to utilize the coherent categorical concept of differential forms on this infinite-dimensional space, which is not possible in ordinary differential geometry.

The Chain Homotopy Operator ($\mathsf{K}$), detailed in [PIZ13a], will also be employed in the generalization of the Moment Map in diffeology, further underscoring its significance.

**14. Homotopy Groups, a Tale of Recursion.** The Cartesian closed nature of the category {Diffeology} allows us to reduce the construction of homotopy groups to the definition of the space of connected components of a diffeological space. Let X be a diffeological space.

DEFINITION (HOMOTOPY). *Two points* $x, x'$ *in* X *are* homotopic*, or* connected by a path*, or simply* connected*, if there exists a smooth path* $\gamma$ *in* X *such that* $\gamma(0) = x$ *and* $\gamma(1) = x'$*, where the space of paths is defined by*

$$\mathrm{Paths}(\mathrm{X}) = \mathscr{C}^\infty(\mathbf{R}, \mathrm{X}),$$

*equipped with the functional diffeology. The endpoints of a path are determined by the maps* $\hat{t} : \mathrm{Paths}(\mathrm{X}) \to \mathrm{X}$ *defined by* $\hat{t}(\gamma) = \gamma(t)$*. The endpoints of a path are the pair* $\mathrm{ends}(\gamma) = (\hat{0}(\gamma), \hat{1}(\gamma))$.

THEOREM (CONNECTED COMPONENTS). *The relation of being homotopic is an equivalence relation on* X*. The equivalence classes of this relation are called* connected components*, and the set of connected components is denoted by*

$$\pi_0(\mathrm{X}) = \mathrm{X}/\mathrm{Homotopy}.$$

*The partition of* X *into connected components is the finest partition that makes* X *the disjoint union of its parts. Therefore, equipped with the quotient diffeology,* $\pi_0(\mathrm{X})$ *is discrete.*

*Moreover, the connected components are also the components for the D-topology, this is why there is no ambiguity to call them simply:* components.

*Proof.* The proof of this theorem uses the classical *concatenation* of paths, defined when a path $\gamma$ ends at the beginning of a path $\gamma'$. However, since this operation does not necessarily preserve smoothness, the paths are modified to be *stationary* at their endpoints by composition with a smooth function that flattens the path at its endpoints (a "bump function"). Because the subspace of stationary paths is homotopically equivalent to the space of paths (via a suitable deformation retraction), this procedure introduces no inconsistency. We refer to [PIZ13, §5.7] for the proof that the connected components are components for the D-topology. ►

The *fundamental group* of a diffeological space X is defined using the fact that Paths(X) is a diffeological space with the functional diffeology, and this applies to all its subspaces. Thus, their sets of connected components are defined without any further specification.

DEFINITION (FUNDAMENTAL GROUP). *The fundamental group of a diffeological space* X *based at* $x \in \mathrm{X}$ *is the space of connected components of the space of loops based at* $x$, *denoted by* $\pi_1(\mathrm{X}, x) = \pi_0(\mathrm{Loops}(\mathrm{X}, x))$, *where*

$$\mathrm{Loops}(\mathrm{X}, x) = \{\ell \in \mathrm{Paths}(\mathrm{X}) \mid \ell(0) = \ell(1) = x\}.$$

*Equipped with the operation* $\mathrm{comp}(\ell) \cdot \mathrm{comp}(\ell') = \mathrm{comp}(\ell \vee \ell')$, *where* comp *denotes the projection to the space of connected components, and* $\ell \vee \ell'$ *denotes the concatenation of loops* $\ell$ *and* $\ell'$ *(with suitable reparametrization to ensure smoothness),* $\pi_1(\mathrm{X}, x)$ *is a group with identity the component of the constant loop at* $x$, *denoted by* $\hat{x} : t \mapsto x$.

This construction of homotopy warrants several remarks:

(1) This definition coincides with the classical definition when X is a manifold.
(2) It also yields the homotopy of the irrational torus: $\pi_0(\mathrm{T}_\alpha) = \{\mathrm{T}_\alpha\}$, $\pi_1(\mathrm{T}_\alpha) = \mathbf{Z}^2$, and $\pi_0(\mathrm{Diff}(\mathrm{T}_\alpha)) = \{\pm 1\} \times \mathbf{Z}$ or $\pi_0(\mathrm{Diff}(\mathrm{T}_\alpha)) = \{\pm 1\}$, depending on whether $\alpha$ is quadratic or not.
(3) While the fundamental group of a diffeological space captures some geometric aspects, it does not fully reflect the finer structure of diffeology. For example, half-lines $\Delta_n$ are contractible and have a trivial homotopy, which obscures their diffeological structures.
(4) As another example, the orbifold $\mathcal{Q}_m = \mathbf{C}/\mathcal{U}_m$ is contractible and thus also has a trivial homotopy. Its detailed diffeological structure is encoded in the structure groupoid, which is not apparent at this level of homotopy.
(5) Nonetheless, homotopy remains a valuable tool in diffeology and deserves further development, as initiated early in the theory [PIZ85].

Now, let us introduce the *higher homotopy groups* recursively.

DEFINITION (HIGHER HOMOTOPY GROUPS). *For any diffeological space* X *and* $x \in \mathrm{X}$, *the* higher homotopy groups *are defined recursively by:*

$$\textit{For all } k \geq 1, \quad \pi_k(\mathrm{X}, x) = \pi_{k-1}(\mathrm{Loops}(\mathrm{X}, x), \hat{x}).$$

*An equivalent presentation of the higher homotopy groups, starting with* $\mathrm{X}_0 = \mathrm{X}$ *and* $x_0 = x$, *uses the iteration:*

$$\mathrm{X}_k = \mathrm{Loops}(\mathrm{X}_{k-1}, x_{k-1}) \quad \textit{with} \quad x_k = [t \mapsto x_{k-1}].$$

*Then, for all* $k > 1$*:*

$$\pi_k(\mathrm{X}, x) = \pi_1(\mathrm{X}_{k-1}, x_{k-1}).$$

This presentation immediately reveals the group structure of the higher homotopy groups. For example: $\pi_2(\mathrm{X}, x) = \pi_1(\mathrm{X}_1, x_1)$, $\mathrm{X}_1 = \mathrm{Loops}(\mathrm{X}, x)$ and $x_1 = \hat{x}$, then $\pi_2(\mathrm{X}, x) = \pi_0(\mathrm{Loops}(\mathrm{X}_1, x_1), \hat{x}_1) = \pi_0(\mathrm{Loops}(\mathrm{Loops}(\mathrm{X}, x), \hat{x}), [t \mapsto [s \mapsto x]])$.

**15. De Rham Cohomology.** Differential calculus, particularly the manipulation of differential forms, provides an intrinsic way to describe geometric objects and physical quantities in differential geometry. Unlike vector fields or other contravariant tensorial objects, differential forms transform naturally under smooth maps. This makes them ideal for studying geometry and physics not only on manifolds but also, as we have seen, on general diffeological spaces.

Differential forms are fundamental in physics: Maxwell's equations in electromagnetism, the description of spacetime curvature in general relativity, the representation of gauge fields in gauge theory, and the modern formulation of classical mechanics with symplectic forms all rely heavily on them. Although differential forms are defined by local properties, their de Rham cohomology groups capture aspects of the global structure. For example, if the first de Rham cohomology group is non-trivial, then the space is not simply connected.

Consider a diffeological space X and its space of differential $k$-forms $\Omega^k(X)$. The vector space of *closed k-forms*, where $d\alpha = 0$, is denoted by $\mathbf{Z}^k_{dR}(X)$. The vector space of *exact k-forms*, where $\alpha = d\beta$, is denoted by $\mathbf{B}^k_{dR}(X)$. The $k$-th *de Rham cohomology group* (which is a vector space) is then defined as usual by:

$$\mathbf{H}^k_{dR}(X) = \mathbf{Z}^k_{dR}(X)/\mathbf{B}^k_{dR}(X).$$

The cohomology of the irrational torus can be computed, yielding $\mathbf{H}^1_{dR}(T_\alpha) = \mathbf{R}$. This is not surprising, as $T_\alpha$ has dimension 1. However, its homotopy is $\mathbf{Z}^2$, which contrasts with what occurs in the differential geometry of manifolds. This discrepancy relates to the non-equivalence of various cohomology theories in diffeology [PIZ24].

**16. The Chain Homotopy Operator.** The *chain-homotopy operator* is a fundamental construction in Cartan-de Rham calculus within diffeology [PIZ13a]. Since it involves differential forms on infinite-dimensional spaces of paths, it extends beyond the scope of classical differential geometry. This smooth linear integro-differential operator connects the de Rham complex of a space to that of its space of paths:

$$\mathbf{K} : \Omega^p(X) \to \Omega^{p-1}(\mathrm{Paths}(X)),$$

satisfying the identity

$$d \circ \mathbf{K} + \mathbf{K} \circ d = \hat{1}^* - \hat{0}^*.$$

The significance of this operator is evident, in particular, in how it streamlines proofs in classical differential geometry, as illustrated by the homotopy invariance of de Rham cohomology, discussed in the next section.

It is also crucial in the construction of the moment map in diffeology, providing a shortcut for the non-exact parasymplectic case and extending this fundamental tool of symplectic geometry beyond manifolds.

The construction of the chain-homotopy operator is based on two operations: the contraction of differential forms by slidings and the diffeological version of the Lie derivative.

(1) *Contraction of differential forms by slidings.* Let $\alpha$ be a $k$-form on X and $\tau$ a germ of a path in Diff(X), centered at the identity, $\tau(0) = \mathbf{1}_X$, called a *sliding*. The contraction of $\alpha$ by $\tau$ is defined, for all plots P in X, by:

$$\iota_\tau(\alpha)(P)_r(v_1,\dots,v_{k-1}) = \alpha\left(\begin{pmatrix} t \\ r \end{pmatrix} \mapsto \tau(t)\circ P(r)\right)_{\binom{0}{r}} \begin{pmatrix} 1 \\ 0 \end{pmatrix}\begin{pmatrix} 0 \\ v_1 \end{pmatrix}\cdots\begin{pmatrix} 0 \\ v_{k-1} \end{pmatrix}$$

The tangent vector $(1,0)$ corresponds to the $\partial/\partial t$ direction, and $(0,v_i)$ corresponds to the tangent vector associated with $v_i$ in $T_r(\mathrm{dom}(P))$.

(2) *Lie derivative by slidings.* Let $\tau$ be a sliding of X and $\alpha$ be a $k$-form. The *Lie derivative* of $\alpha$ by $\tau$ is defined by:

$$\mathfrak{L}_\tau(\alpha) = \left.\frac{\partial\tau(t)^*(\alpha)}{\partial t}\right|_{t=0}.$$

(3) *Integration of forms on paths.* Consider a $k$-form $\alpha$ on X, and the operator $\hat{t} : \mathrm{Paths}(X) \to X$ defined by $\hat{t}(\gamma) = \gamma(t)$. Let $\Phi$ be the integral operator from $\Omega^k(X)$ to $\Omega^k(\mathrm{Paths}(X))$, defined, for all plots P in Paths(X), by:

$$\Phi(\alpha)(P)_r(v_1,\dots,v_k) = \int_0^1 \hat{t}^*(\alpha)(P)_r(v_1,\dots,v_k)\,dt.$$

Then, the integro-differential operator $\mathsf{K}$ is defined as the composition of $\iota_\tau$ and $\Phi$:

$$\mathsf{K} = \iota_\tau\circ\Phi, \quad \mathsf{K}\in L^\infty(\Omega^k(X),\Omega^{k-1}(\mathrm{Paths}(X))).$$

The Chain Homotopy identity is derived by applying the Lie derivative to the integral operator $\Phi$ with respect to the 1-parameter group $\tau$, and then applying the generalized Cartan formula (op. cit.):

$$\underbrace{\mathfrak{L}_\tau\Phi(\alpha)}_{\hat{1}^*(\alpha)-\hat{0}^*(\alpha)} = d[\iota_\tau(\Phi(\alpha))] + \underbrace{\iota_\tau(d(\Phi(\alpha)))}_{d\circ\Phi=\Phi\circ d} = \underbrace{d[\iota_\tau(\Phi(\alpha))]+\iota_\tau[\Phi(d\alpha)]}_{d[\mathsf{K}(\alpha)]+\mathsf{K}[d\alpha]}$$

This construction, based on the Cartan-Lie derivative, demonstrates how diffeology enables us to address problems directly, without involving unnecessary objects, particularly vector fields. Their contravariant nature leads to multiple, non-equivalent definitions within diffeology. Their inclusion should be avoided whenever possible, as is the case here, to maintain the general applicability of the theorems. Avoiding vector fields and Lie algebras in diffeology is not merely a stylistic choice but a requirement for our constructions and theorems to be applicable throughout the entire category of diffeological spaces.

**17. Application to the Homotopy Invariance of de Rham Cohomology.** Perhaps the most compelling application of the chain-homotopy operator lies in the proof of the homotopy invariance of de Rham cohomology. The use of infinite-dimensional spaces significantly simplifies this proof, which is notoriously intricate and tedious when approached using standard differential geometry techniques.

THEOREM (HOMOTOPY INVARIANCE OF DE RHAM COHOMOLOGY). *Let* X *and* X′ *be two diffeological spaces, and let* $t\mapsto f_t$ *be a smooth path in* $\mathscr{C}^\infty(X',X)$. *For any closed differential form* $\alpha\in\Omega^k(X)$, *the forms* $f_1^*(\alpha)$ *and* $f_0^*(\alpha)$ *are cohomologous.*

*Proof.* Let $\mathrm{F} : \mathrm{X}' \to \mathrm{Paths}(\mathrm{X})$ be the map defined by $\mathrm{F}(x') = [t \mapsto f_t(x')]$. Applying $\mathrm{F}^*$ to the chain-homotopy operator identity, we obtain:

$$\underbrace{\mathrm{F}^*(d\mathbf{K}\alpha}_{d(\mathrm{F}^*(\mathbf{K}\alpha)=\beta)} \quad + \quad \mathbf{K}\cancel{d\alpha}) \quad = \quad \underbrace{\mathrm{F}^*(\hat{1}^*(\alpha))}_{(\hat{1}\circ \mathrm{F})^*(\alpha) = f_1^*(\alpha)} \quad - \quad \underbrace{\mathrm{F}^*(\hat{0}^*(\alpha))}_{(\hat{0}\circ \mathrm{F})^*(\alpha) = f_0^*(\alpha)}$$

That is, $f_1^*(\alpha) = f_0^*(\alpha) + d\beta$. ▶

This example demonstrates the utility of diffeology, even within the context of traditional differential geometry, as it provides a rigorous proof applicable when X is a manifold. Furthermore, the theorem's applicability extends to all diffeological spaces, a point worth emphasizing.

## 5. The Moment Maps in Diffeology

The moment map is a fundamental tool in the study of symmetries in symplectic geometry. Introduced by J.-M. Souriau in the 1960s, its significance spans a wide range of applications, from the Noether-Souriau theorem on conserved quantities in classical mechanics to geometric quantization and the classification of elementary particles as homogeneous symplectic spaces.

However, the application of moment maps becomes significantly more challenging in infinite-dimensional settings, such as those encountered in models of field theory. For nearly half a century, mathematicians and physicists have developed such models, equipped with candidate symplectic structures and infinite-dimensional groups of symmetries. Notable examples of this work, where heuristic approaches to moment maps have been proposed, include the work of M. Atiyah [AB82] and S. Donaldson [Dnl99].

Diffeology provides a natural framework to successfully address these challenges and generalize the concept of the moment map in a unified and rigorous manner. By offering a consistent approach, diffeology obviates the need for ad hoc constructions and heuristics. Moreover, it applies effectively to problems involving both infinite-dimensional and singular spaces [PIZ10]. This framework allows us, in particular, to work directly with the *space of momenta*, a key element in this generalization.

Recall that a *diffeological group* is a diffeological space equipped with a group law such that the product and the inverse maps are smooth. For example, every group of diffeomorphisms Diff(X) of a diffeological space X is a diffeological group when equipped with the functional diffeology.

Definition (Momentum, Momenta). *A* momentum *of a diffeological group* G *is a left-invariant* 1*-form. The set of momenta forms a diffeological vector space under the functional diffeology. We denote it by:*

$$\mathscr{G}^* = \{\varepsilon \in \Omega^1(\mathrm{G}) \mid \forall g \in \mathrm{G}, \mathrm{L}(g)^*(\varepsilon) = \varepsilon\},$$

*where* L *denotes left multiplication:* $\mathrm{L}(g) \colon g' \mapsto gg'$.

It is important to emphasize that the space of momenta $\mathscr{G}^*$ is *not* defined by duality with a presumed Lie algebra, even if the notation might suggest it. The introduction of Lie algebras is unnecessary in the definition of the moment map, as we shall see, and indeed obscures the fundamental simplicity of this object.

**18. The Simplest Case.** Recall that a *parasymplectic form* on a diffeological space X is any closed 2-form on X, without regard to its regularity. Let G be a diffeological group acting smoothly on X and preserving $\omega$, i.e., $\underline{g}^*(\omega) = \omega$ for all $g \in \mathrm{G}$.

Assume that $\omega$ is exact, possessing an invariant primitive $\alpha$: $\omega = d\alpha$ and $\underline{g}^*\alpha = \alpha$ for all $g \in \mathrm{G}$.

DEFINITION (SIMPLE MOMENT MAP). *The map* $\mu \colon \mathrm{X} \to \mathcal{G}^*$ *defined by*

$$\mu(x) = \hat{x}^*(\alpha),$$

*where* $\hat{x} : \mathrm{G} \to \mathrm{X}$ *is the* orbit map *given by* $\hat{x}(g) = \underline{g}(x)$, *is smooth, i.e.,* $\mu \in \mathcal{C}^\infty(\mathrm{X}, \mathcal{G}^*)$. *This map* $\mu$ *is a* moment map *for* $\omega$.

Note that this is one moment map for $\omega$, associated with the chosen invariant primitive $\alpha$.

**19. The General Case.** Closed invariant 2-forms are not always exact, and even if exact, may not possess an invariant primitive. Diffeology addresses this issue directly by employing the chain-homotopy operator to reduce to the simplest case. Let $(\mathrm{X}, \omega)$ be a connected parasymplectic space.

(1) *The Paths Moment Map.* To reduce the general case to the simple one, we shift our perspective from the space X to the space of paths, Paths(X). Let us define a 1-form $\lambda$ on Paths(X) by applying the chain-homotopy operator to $\omega$:

$$\lambda := \mathsf{K}\omega.$$

By the fundamental property of the chain-homotopy operator, the exterior derivative of this 1-form is $\varpi = d\lambda = \hat{1}^*(\omega) - \hat{0}^*(\omega)$. It is established that if G preserves $\omega$, then its action on paths preserves $\lambda$ (see [PIZ13, §9.3]). We are therefore precisely in the "simplest case" discussed above, but on the space of paths: we have a group G acting on Paths(X) and preserving a 1-form $\lambda$. We can thus directly apply the simple definition of the moment map.

DEFINITION (THE PATHS MOMENT MAP). *The* paths moment map *is the moment map for the action of* G *on* $(\mathrm{Paths}(\mathrm{X}), d\lambda)$ *associated with the invariant primitive* $\lambda$*:*

$$\Psi : \mathrm{Paths}(\mathrm{X}) \to \mathcal{G}^* \quad \textit{with} \quad \Psi(\gamma) = \hat{\gamma}^*(\lambda).$$

The paths moment map satisfies two fundamental properties:

(1) $\Psi$ is equivariant with respect to the coadjoint action of G:[17]

$$\Psi \circ g_* = \mathrm{Ad}_*(g) \circ \Psi \ \text{ where } \ g_*(\gamma) = \underline{g} \circ \gamma.$$

(2) $\Psi$ is additive: for composable paths $\gamma$ and $\gamma'$,

$$\Psi(\gamma \vee \gamma') = \Psi(\gamma) + \Psi(\gamma').$$

[17]Recall that, for all $g, k \in \mathrm{G}$, $\mathrm{ad}(g)(k) = gkg^{-1}$, and for all $\alpha \in \mathcal{G}^*$, $\mathrm{Ad}_*(g)(\alpha) = [\mathrm{ad}(g)]_*(\alpha)$.

(2) *The Two-Points Moment Map.* The *two-points moment map* is the projection of the paths moment map onto $X \times X$. This projection provides a first step towards defining the moment map of $\omega$ on the space $X$:

$$\psi : X \times X \to \mathscr{G}^*/\Gamma, \quad \text{with} \quad \psi(x, x') = \text{class}(\Psi(\gamma)),$$

where class : $\mathscr{G}^* \to \mathscr{G}^*/\Gamma$ and for any path $\gamma$ such that $\text{ends}(\gamma) = (x, x')$. The subspace $\Gamma \subset \mathscr{G}^*$ is the obstruction to the projection ends, i.e.,

$$\Gamma = \{\Psi(\ell) \mid \ell \in \text{Loops}(X)\}.$$

It is a subgroup under addition, consisting of closed $\text{Ad}_*$-invariant momenta, and $\mathscr{G}^*/\Gamma$ is regarded as an Abelian diffeological group. It is the *holonomy* of the action of G on $(X, \omega)$, representing the obstruction to the action of G being *Hamiltonian.* The action of G on $(X, \omega)$ is Hamiltonian precisely when $\Gamma = \{0\}$ and $\psi$ takes values in $\mathscr{G}^*$.

The 2-points moment map $\psi$ remains G-equivariant and satisfies an additive property, making it a *Chasles cocycle*:

$$\psi(\underline{g}(x), \underline{g}(x')) = \text{Ad}_*(\psi(x, x')) \;\text{ and }\; \psi(x, x') + \psi(x', x'') = \psi(x, x'').$$

(3) *The One-Point Moment Map.* A *one-point moment map* is a solution $\mu$ to the equation

$$\psi(x, x') = \mu(x') - \mu(x) \quad \text{with} \quad \mu : X \to \mathscr{G}^*/\Gamma.$$

Since X is connected, the solutions differ by a constant and are given by:

$$\mu(x) = \psi(x_0, x) + c, \;\text{ for some choice of } x_0 \in X \text{ and } c \in \mathscr{G}^*/\Gamma.$$

The moment map $\mu$ is no longer $\text{Ad}_*$-equivariant but $\theta$-affine $\text{Ad}_*$-equivariant:

$$\mu(\underline{g}(x)) = \text{Ad}_*(\mu(x)) + \theta(g), \quad \text{where} \quad \theta(g) = \psi(x_0, \underline{g}(x_0)) - \Delta(c)(g).$$

The map $\theta : G \to \mathscr{G}^*/\Gamma$ is a 1-cocycle of G with values in $\mathscr{G}^*/\Gamma$, twisted by the coadjoint action. We denote $\theta \in H^1(G, \mathscr{G}^*/\Gamma)$. The map $\Delta c$ is a coboundary for this cohomology, $\Delta(c)(g) = \text{Ad}_*(g)(c) - c$. Denoting by $\text{Ad}_*^\theta : g \mapsto \text{Ad}_*(g) + \theta(g)$, the moment map $\mu$ is $\text{Ad}_*^\theta$-equivariant.

This construction warrants a couple of remarks.

(1) When X is a manifold and G a Lie group, if the action is Hamiltonian, then the 1-point moment map coincides with Souriau's definition in [Sou70], and the cocycle coincides with Souriau's cocycle. Hence, we retain this terminology. In contrast to Souriau's construction, we avoid two elements: the Lie algebra of G and the need to solve any differential equations. This is noteworthy and underscores that diffeology adopts a global perspective.

(2) This construction applies, in particular, to the entire group of automorphisms $G_\omega = \text{Diff}(X, \omega)$. The momenta and algebraic objects arising in this construction are indexed by $\omega$ and are termed *universal*: $\Psi_\omega$, $\psi_\omega$, $\mu_\omega$, $\Gamma_\omega$, and $\theta_\omega$. They encapsulate all information about the symmetries of $\omega$ encoded in the diffeology of X.

**20. Examples: Atiyah, Donaldson and Virasoro.** To illustrate the construction of the moment map in diffeology, we will now present three key examples of moment maps in infinite-dimensional settings. Further examples, including those in finite dimensions and on singular spaces, can be found in [PIZ10].

(1) *The Moment of Imprimitivity.* Let M be a manifold. The Abelian group $\mathrm{A} = \mathscr{C}^\infty(\mathrm{M}, \mathbf{R})$ acts on the cotangent bundle $\mathrm{T}^*\mathrm{M}$ by vertical translation: $\underline{f}(q,p) = (q, p + df(q))$ for $f \in \mathrm{A}$. This action preserves the canonical symplectic form $\omega = dp \wedge dq$. For each point $q \in \mathrm{M}$, let us define the *evaluation function* $\delta_q \in \mathscr{C}^\infty(\mathrm{A}, \mathbf{R})$ by $\delta_q(f) = f(q)$. This function is smooth and not invariant under the group action on A, but its exterior differential, $d\delta_q$, is a left-invariant 1-form on A, i.e., an element of the space of momenta $\mathscr{A}^*$. The moment map for this action is then given by:

$$\mu : (q,p) \mapsto d\delta_q \in \mathscr{A}^*.$$

Note that $\delta_q$ can be identified with the Dirac distribution centered at $q$. The moment map is thus the differential of a distribution, a frequent occurrence in this generalized setting.

(2) *Topological Field Theory Example.* Let $\Sigma$ be an oriented surface. The Abelian group $\mathrm{A} = \mathscr{C}^\infty(\Sigma, \mathbf{R})$ acts on the space of 1-forms $\Omega^1(\Sigma)$, preserving the 2-form $\omega(\alpha,\beta) = \int_\Sigma \alpha \wedge \beta$. For all $f \in \mathrm{A}$ and $\alpha \in \Omega^1(\Sigma)$, let $\underline{f}(\alpha) = \alpha + df$. This example is drawn from the literature [AB82, Dnl99]. Its moment map in diffeology is:

$$\mu : \alpha \mapsto d\left[f \mapsto \int_\Sigma f \, d\alpha\right] \in \mathscr{A}^*.$$

Here again, the moment map is the differential of a distribution: $[d\alpha] : f \mapsto \int_\Sigma f \, d\alpha$. Heuristically, $d\alpha$ is often considered as the moment map, but diffeology precisely defines its status.

(3) *Virasoro.* Consider the space $\mathrm{X} = \mathrm{Imm}(\mathrm{S}^1, \mathbf{R}^2)$ of immersions of $\mathrm{S}^1 = \mathbf{R}/2\pi\mathbf{Z}$ in $\mathbf{R}^2$, and $\omega = d\alpha$, the exterior derivative of

$$\alpha(\delta x) = \int_0^{2\pi} \frac{1}{\|\dot{x}(t)\|^2} \langle \ddot{x}(t) | \delta\dot{x}(t) \rangle \, dt,$$

where $x \in \mathrm{Imm}(\mathrm{S}^1, \mathbf{R}^2)$ and $\delta x$ is a variation of $x$ induced by a parameter $s \mapsto x_s$, i.e., $x_0 = x$ and $\delta x = \left.\frac{dx_s}{ds}\right|_{s=0}$. The group $\mathrm{Diff}^+(\mathrm{S}^1)$ acts on X by reparametrization: $\underline{\phi}(x) = x \circ \phi^{-1}$.

On the connected component of the standard immersion $t \mapsto (\cos(t), \sin(t))$, the moment map is, up to a constant:

$$\mu(x)(\mathrm{P})_r(\delta r) = \int_0^{2\pi} \left\{ \frac{\|x''(u)\|^2}{\|x'(u)\|^2} - \frac{d^2}{du^2} \log \|x'(u)\|^2 \right\} \delta u \, du,$$

where $\mathrm{P} : \mathrm{dom}(\mathrm{P}) \to \mathrm{Diff}_+(\mathrm{S}^1)$ is an $n$-plot, $r \in \mathrm{dom}(\mathrm{P})$, $\delta r \in \mathbf{R}^n$, $u = \phi^{-1}(t)$, $t$ is the parameter of $x \in \mathrm{Imm}(\mathrm{S}^1, \mathbf{R}^2)$, and $\delta u = \mathrm{D}(r \mapsto u)(r)(\delta r)$.

The moment map is not $\mathrm{Ad}_*$-equivariant and defines a non-trivial cocycle of $\mathrm{Diff}_+(S^1)$ cohomologous to $\theta$, defined by:

$$\theta(g)(P)_r(\delta r) = \int_0^{2\pi} \frac{3\gamma''(u)^2 - 2\gamma'''(u)\gamma'(u)}{\gamma'(u)^2}\delta u \, du,$$

where $g \in \mathrm{Diff}^+(S^1)$ and $\gamma = g^{-1}$, revealing the integrand of the right-hand side as the celebrated *Schwarzian derivative*. This result provides a powerful illustration of the diffeological approach: it identifies the geometric origin of the Schwarzian derivative — a classical invariant under Möbius transformations — as precisely the cocycle measuring the lack of equivariance of the moment map in this natural, infinite-dimensional symplectic construction.[18] Furthermore, it is this very cocycle that, through the lens of the prequantum construction in diffeology (as we will see in the next section), gives rise to the famous Virasoro group — a central extension of $\mathrm{Diff}_+(S^1)$ — and its associated Lie algebra, the Virasoro algebra.

To conclude this chapter on generalized moment maps, it is evident that diffeology provides a simple, rigorous, and effective framework for dealing with complex situations involving, as we have seen, infinite-dimensional spaces and distributions, which pose significant challenges for ordinary differential geometry. While functional analysis offers some tools, diffeology provides a more direct and natural approach, yielding clear and expected results. These examples strongly demonstrate the power and utility of diffeology.

## 6. Prequantum Groupoid and Bundles

Geometric quantization is a fundamental program in mathematical physics aiming to bridge classical and quantum mechanics. A crucial first step is prequantization, which, for a classical system described by a symplectic manifold $(M, \omega)$, involves constructing a principal $U(1)$-bundle over M with curvature $\omega$. Applying this framework to more general spaces, such as those with singularities or infinite-dimensional phase spaces encountered in field theories, presents significant challenges for traditional manifold-based methods. Diffeology provides a powerful and natural setting to generalize this construction. By leveraging the capabilities of functional diffeology, we can endow the space of paths $\mathrm{Paths}(X)$ with a smooth structure and perform calculus on it via the chain-homotopy operator. This allows us to define a prequantum groupoid $(\mathbf{T}_\omega, \boldsymbol{\lambda})$ as a diffeological quotient of $\mathrm{Paths}(X)$, whose structure intrinsically encodes the classical system $(X, \omega)$ and its quantum phase information (via the isotropy group isomorphic to the torus of periods). This construction, detailed in Article 22 (see also [PIZ25b]), serves as a prime example of how diffeology enables the extension of key geometric and physical concepts to spaces far beyond the reach of classical differential geometry.

**21. Diffeological Fiber Bundles.** The question of how to define fiber bundles in diffeology arose in 1983 with the computation of the homotopy of the irrational torus [DIZ83]. In the case of the projection $\pi : T^2 \to T_\alpha$, it might appear that we have a fiber bundle with fiber $\mathbf{R}$. If this were the case, we could apply the exact long homotopy sequence, and since $\mathbf{R}$ is contractible, we would immediately obtain $\pi_1(T_\alpha) = \mathbf{Z}^2$ (and $\pi_k(T_\alpha) = 0$ for all $k > 1$). However, this approach fails because $T_\alpha$ is topologically trivial (i.e., it inherits the trivial D-topology by quotient) and the projection $\pi$ is not locally trivial.

[18] A formal construction in the rigorous diffeology framework, not an had hoc heuristic addendum.

$$\begin{array}{rcl} & T^2 & \leftarrow \text{Total space} \\ \text{Fiber } \mathbf{R} \rightarrow & \Big\downarrow \pi & \\ & T_\alpha & \leftarrow \text{Base space} \end{array}$$

It was therefore necessary to revise the classical concept of fiber bundle to adapt it to the requirements of diffeology:

(1) The projection of a diffeological group onto its quotient by any subgroup is a diffeological fibration.
(2) Diffeological fiber bundles satisfy the exact long homotopy sequence.

This problem was resolved in 1985 by replacing local triviality over the base with *local triviality along plots* [PIZ85].

DEFINITION (FIBER BUNDLES). *The projection* $\pi : Y \to X$ *is a (diffeological)* fibration,[19] *with fiber* F*, if for every plot* $P : U \to X$*, the pullback*

$$\mathrm{pr}_1 : P^*(Y) \to U \quad \textit{where} \quad P^*(Y) = \{(r, y) \in U \times Y \mid P(r) = \pi(y)\},$$

*is locally trivial with fiber* F.

That is, for every point in the domain U, there exists an open neighborhood V such that there is a diffeomorphism $\varphi : V \times F \to \mathrm{pr}_1^{-1}(V) \subset P^*(Y)$, satisfying $\mathrm{pr}_1 \circ \varphi = \mathrm{pr}_1$. This is summarized in the following diagram:

$$\begin{array}{ccccc} V \times F & \xrightarrow{\ \varphi\ } & P^*(Y) \restriction V & \xrightarrow{\ \mathrm{pr}_2\ } & Y \\ & \searrow^{\mathrm{pr}_1} & \Big\downarrow \mathrm{pr}_1 & & \Big\downarrow \pi \\ & & V & \xrightarrow[P \restriction V]{} & X \end{array}$$

This concept of local triviality along plots represents a departure from the classical definition of fiber bundles, where triviality is required over open sets of the base space itself. In diffeology, this weaker condition allows for a more flexible treatment of fiber bundles on spaces with non-standard smooth structures, and ensures that fiber bundles satisfy the following:

1. When X is a manifold, local triviality along plots is equivalent to local triviality. If the fiber F is also a manifold, then the total space Y is a manifold, and $\pi$ is a fibration in the classical sense. Hence, diffeological fiber bundles fully and faithfully extend their counterparts in ordinary differential geometry.

[19] The terms “fibration” and “fiber bundle” are used interchangeably.

2. The projection $\pi\colon G \to G/H$, where G is a diffeological group and $H \subset G$ is any subgroup,[20] is a diffeological fiber bundle. It satisfies local triviality along plots. This is the case for the projection $\pi : T^2 \to T_\alpha = T^2/\mathscr{S}_\alpha$.
3. Diffeological fibrations $\pi : Y \to X$ with fiber F satisfy the following exact long homotopy sequence, where $x \in X$ and $y \in \pi^{-1}(x)$:

$$\cdots \to \pi_n(F, y) \to \pi_n(Y, y) \to \pi_n(X, x) \to \pi_{n-1}(F, y) \to \cdots$$
$$\cdots \to \pi_0(F, y) \to \pi_0(Y, y) \to \pi_0(X, x) \to 0.$$

This property is crucial for preserving essential geometric tools and intuition within diffeology.

Thus, this definition generalizes the classical concept of fiber bundles, providing the necessary flexibility within diffeology without contradicting classical theory in the familiar setting of manifolds. Covering spaces, an important class of fiber bundles, also extend naturally to diffeology.

DEFINITION (COVERING). *A* covering space *of a connected diffeological space* X *is a diffeological fiber bundle* $\pi : \hat{X} \to X$*, where the fiber* $\pi^{-1}(x)$ *over any point* $x$ *in* X *is a discrete diffeological space.*

Covering spaces are classified by the following theorem:

THEOREM (UNIVERSAL COVERING). *Every connected diffeological space* X *admits a simply connected universal covering* $\pi : \tilde{X} \to X$*, unique up to isomorphism, with fiber the fundamental group* $\pi_1(X, x)$*. Every other connected covering is obtained as a quotient of* $\tilde{X}$ *by a subgroup of* $\pi_1(X, x)$*.*

The construction of $\tilde{X}$ is sufficiently interesting to warrant a sketch; further details can be found in [PIZ85, PIZ13].

*Sketch of Proof.* Let $\mathscr{X}$ be the quotient space

$$\mathscr{X} = \text{Paths}(X)/\text{fixed-ends homotopy},$$

where the fixed-ends homotopy relation is defined as follows: paths $\gamma$ and $\gamma'$ are *fixed-ends homotopic* if:

(1) They have the same endpoints: $\text{ends}(\gamma) = \text{ends}(\gamma')$.
(2) There exists a path $[s \mapsto \gamma_s]$ in Paths(X) such that
   (a) $\text{ends}(\gamma_s) = \text{ends}(\gamma) = \text{ends}(\gamma')$ for all $s$,
   (b) $\gamma_0 = \gamma$ and $\gamma_1 = \gamma'$.

As a quotient of a diffeological space, $\mathscr{X}$ is a diffeological space. It carries the structure of a *diffeological groupoid* (op. cit.) whose objects are the points in X and whose morphisms are the fixed-ends homotopy classes:

$$\text{Obj}(\mathscr{X}) = X \quad \text{and} \quad \text{Mor}_{\mathscr{X}}(x, x') = \text{Paths}(X, x, x')/\text{fixed-ends homotopy}.$$

The isotropy group at each point $x$ is the fundamental group at that point:

$$\text{Mor}_{\mathscr{X}}(x, x) = \pi_1(X, x).$$

[20]Any subgroup equipped with the subset diffeology is a diffeological group.

Fixing a base point $x$ in X, let

$$\tilde{\mathrm{X}} = \mathrm{Mor}_{\mathscr{X}}(x, *) = \mathrm{Paths}(\mathrm{X}, x)/\text{fixed-ends homotopy}.$$

This space is equipped with the subset diffeology induced by the quotient diffeology on the space of paths. Let

$$\pi : \tilde{\mathrm{X}} \to \mathrm{X} \quad \text{where} \quad \pi(\mathrm{class}(\gamma)) = \hat{1}(\gamma).$$

The space $\tilde{\mathrm{X}}$ is simply connected, the projection $\pi$ is a fibration, and the action of the isotropy group $\pi_1(\mathrm{X}, x) = \mathrm{Mor}_{\mathscr{X}}(x, x)$ makes $\tilde{\mathrm{X}}$ a principal bundle.[21] ►

The universality of $\tilde{\mathrm{X}}$, i.e., the fact that any other connected covering is a quotient of the universal covering by a subgroup of the fundamental group, follows from the *monodromy theorem*.

THEOREM (MONODROMY). *Let $\pi : \tilde{\mathrm{X}} \to \mathrm{X}$ be the universal covering of a diffeological space* X. *Let* X′ *be a simply connected diffeological space and $f \in \mathscr{C}^\infty(\mathrm{X}', \mathrm{X})$. Choose $x' \in \mathrm{X}'$ and $x = f(x')$. For any $\tilde{x} \in \pi^{-1}(x)$, there exists a unique lifting $\tilde{f} \in \mathscr{C}^\infty(\mathrm{X}', \tilde{\mathrm{X}})$ such that $\pi \circ \tilde{f} = f$ and $\tilde{f}(x') = \tilde{x}$.*

The construction of the universal covering space within diffeology represents a significant advancement in our ability to rigorously define and work with this important concept. While mathematicians often possess a strong intuitive understanding of universal coverings, particularly in familiar settings like manifolds, the challenge lies in translating this intuition into a mathematically sound construction, especially when dealing with infinite-dimensional spaces or spaces with non-standard smooth structures. Classical topological methods frequently fall short in these scenarios, leading to ambiguities or limitations. Diffeology overcomes these limitations by providing a general and rigorous framework for constructing universal coverings, offering a powerful tool for exploring the topology and geometry of a much broader class of spaces. In particular, diffeology enables us to define and construct universal coverings in cases where the base space has a trivial D-topology or non-standard smooth structure, providing a level of generality and rigor often unattainable with classical methods.

Thanks to this construction, we know, for example, that the identity component of every group of diffeomorphisms of every diffeological space possesses a unique simply connected universal covering. Consider, for instance, the group of direct diffeomorphisms of the circle $\mathrm{S}^1 = \mathbf{R}/2\pi\mathbf{Z}$, equipped with its functional diffeology. Its universal covering is the group of diffeomorphisms of the real line $\mathbf{R}$ satisfying $f'(x) > 0$ and $f(x + 2\pi) = f(x) + 2\pi$. Originally an example due to P. Donato for homogeneous spaces [Don84], this result now holds within the more general setting of diffeological spaces. Moreover, this construction applies not only to infinite-dimensional diffeological spaces but also to singular ones, even those with a trivial D-topology, such as the irrational torus $\mathrm{T}_\alpha$, whose universal covering is $\mathbf{R}$ and whose first homotopy group is $\mathbf{Z} + \alpha\mathbf{Z} \subset \mathbf{R}$. It is important to understand that when the $\pi_1$ is discrete in $\tilde{\mathrm{X}}$, it is discrete in the diffeological sense: the plots are locally constant, as illustrated by the case of $\mathrm{T}_\alpha$.

[21] Principal bundles are a crucial class of fiber bundle; essentially, they arise when the fibers are the orbits of a smooth free action of a diffeological group [PIZ13, §8.11].

**22. Prequantum Groupoid & Bundle: A Quotient of Paths.** The concept of a prequantum bundle is fundamental in geometric quantization, providing the initial step in a geometric framework for quantizing classical systems. This framework aims to implement Dirac's program, which seeks to associate a Hilbert space to each symplectic manifold (representing the space of solutions of a dynamical system) and to map the algebra of real functions to the algebra of operators, such that the Poisson bracket is mapped to the Lie bracket. The space of smooth functions on the symplectic manifold can be interpreted as the Lie algebra of the group of automorphisms (symplectomorphisms) via the symplectic gradient. Thus, Dirac's program can be more precisely stated as the search for a unitary (or projective) representation of the group of symplectomorphisms into the group of unitary (or projective) transformations of the Hilbert space. A primary difficulty is addressed by replacing the symplectomorphisms with the automorphisms of a prequantum bundle over the symplectic manifold, known as quantomorphisms. Despite successes, such as in the harmonic oscillator, geometric quantization has limitations, e.g., the hydrogen atom. Nevertheless, research continues, particularly as many physical systems involve spaces beyond smooth manifolds, necessitating a generalized approach to geometric quantization, beginning with a generalized prequantum bundle.

Mathematically, a prequantum bundle over a parasymplectic space $(\mathrm{X}, \omega)$ is analogous to integrating a closed 1-form. Consider a closed 1-form $\alpha$ on a diffeological space X. If X is simply connected, $\alpha$ is exact, possessing a primitive:[22]

$$f(x) = \int_{\gamma} \alpha := \int_0^1 \alpha(\gamma)_t(1)\, dt,$$

where $\gamma$ connects a base point $x_0$ to $x$. Indeed, $\alpha = df$. If $\alpha$ is not exact, consider its pullback $\tilde{\alpha}$ to the universal covering, which is simply connected. Then, $\tilde{\alpha}$ is exact; let $\tilde{f}$ be a primitive. The obstruction for $\tilde{f}$ to descend to X is its values on the orbits of the fundamental group, i.e., the integral of $\alpha$ over loops in X. Let

$$\mathrm{P}_\alpha = \left\{ \int_\ell \alpha \mid \ell \in \mathrm{Loops}(\mathrm{X}, x_0) \right\},$$

denote the *group of periods* of $\alpha$. If $\mathrm{P}_\alpha$ is discrete,[23] then define the *torus of periods* as

$$\mathrm{T}_\alpha = \mathbf{R}/\mathrm{P}_\alpha.$$

THEOREM (INTEGRATION OF CLOSED ■-FORMS). *There exists a smooth map* $f \in \mathscr{C}^\infty(\mathrm{X}, \mathrm{T}_\alpha)$ *such that* $\alpha = f^*(\theta)$, *where* $\theta$ *is the canonical volume form on* $\mathrm{T}_\alpha$, *the pushforward of* $dt$ *on* $\mathbf{R}$. *In fact,* $\tilde{f}$ *factors through a covering* $\hat{\mathrm{X}}_\alpha$, *a quotient of* $\tilde{\mathrm{X}}$ *by the kernel of* $f$, *with* $\tilde{f} = \hat{f} \circ \mathrm{pr}_\alpha$.

[22] The notation $\alpha(\gamma)_t(1)$ applies to general diffeological spaces; for manifolds, it is equivalent to $\alpha(\dot{\gamma}(t))$.

[23] Discreteness in diffeology means plots are locally constant. In $\mathbf{R}$, every proper subgroup is discrete. Thus, $\mathrm{P}_\alpha$ being discrete is equivalent to $\mathrm{P}_\alpha \neq \mathbf{R}$.

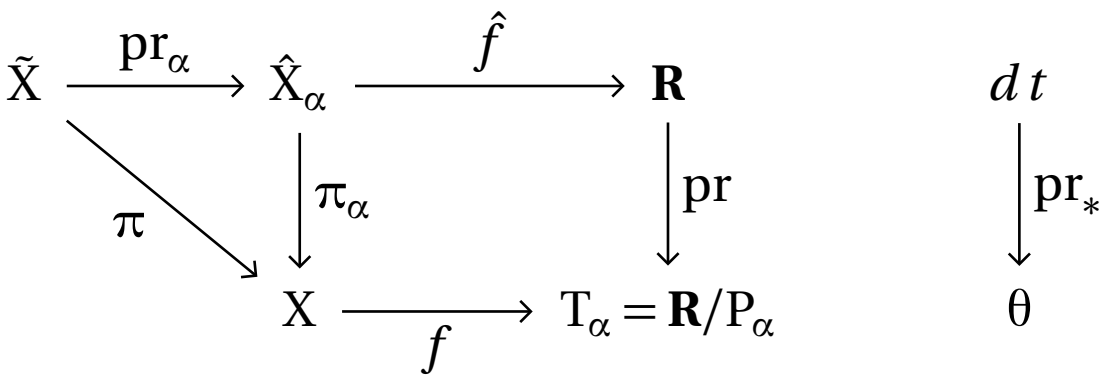

The diffeological fibration $\pi_\alpha : \hat{X}_\alpha \to X$, with the real function $\hat{f}$, integrates the closed 1-form $\alpha$. When $P_\alpha = \{0\}$, we recover the exact case.

When X is a standard manifold, $P_\alpha$ is always discrete, generalizing the familiar case where $P_\alpha = a\mathbf{Z}$. In that instance, $T_\alpha$ is the circle of perimeter $a$. This construction extends to all diffeological spaces, provided $P_\alpha$ is discrete ($P_\alpha \neq \mathbf{R}$).

The prequantum bundle over a parasymplectic space $(X, \omega)$ is analogous, with the covering $\hat{X}_\alpha$ replaced by a principal fiber bundle $\pi : Y \to X$, with structure group $T_\omega$ (the torus of periods of $\omega$), and the real function $f$ replaced by a connection 1-form $\lambda$ with curvature $\omega$.

In what follows, we assume that X is a simply connected diffeological space, i.e., $\pi_0(X) = \{X\}$ and $\pi_1(X, x) = \pi_0(\mathrm{Loops}(X, x)) = \{0\}$ for any base point $x$. Let $\omega$ be a closed 2-form on X. As in the construction of the moment map in diffeology, we employ the chain-homotopy operator $\mathsf{K}$ to transition from X to Paths(X). The 1-form $\mathsf{K}\omega \in \Omega^1(\mathrm{Paths}(X))$, restricted to Loops(X, $x$), is closed because $d\omega = 0$, and on Loops(X, $x$), $\hat{0} = \hat{1}$, implying $\hat{0}^*(\omega \restriction \mathrm{Loops}(X, x)) - \hat{1}^*(\omega \restriction \mathrm{Loops}(X, x)) = 0$. Let $P_\omega$ and $T_\omega$ denote the group of periods and the torus of periods, respectively, of $\mathsf{K}\omega$ on Loops(X, $x$):

$$P_\omega = \left\{ \int_\sigma \mathsf{K}\omega \;\middle|\; \sigma \in \mathrm{Loops}\,(\mathrm{Loops}(X, x), \hat{x}) \right\} \quad \text{and} \quad T_\omega = \mathbf{R}/P_\omega,$$

where $\hat{x}$ represents the constant loop based at $x$. For any pair of paths $\gamma$ and $\gamma'$ that share the same endpoints, i.e., $\gamma(0) = \gamma'(0)$ and $\gamma(1) = \gamma'(1)$, the juxtaposition $\gamma \vee \bar{\gamma}'$, where $\bar{\gamma}'(t) = \gamma'(1-t)$, forms a loop in X based at the origin of the two paths. Define the following equivalence relation on Paths(X, $x$):

$$\gamma \underset{\omega}{\sim} \gamma' \iff \mathrm{ends}(\gamma) = \mathrm{ends}(\gamma') \text{ and } \int_{\hat{x}}^{\gamma \vee \bar{\gamma}'} \mathsf{K}\omega \in P_\omega,$$

where the integral is taken over a path in Loops(X, $x$) connecting $\hat{x}$ to $\gamma \vee \bar{\gamma}'$. This relation is well-defined because every path is fixed-ends homotopic to a stationary path, and the integrals of closed forms are constant on homotopy.[24] Let $\mathscr{Y}$ be the quotient

$$\mathscr{Y} = \mathrm{Paths}(X)/\underset{\omega}{\sim} \text{ with } \mathrm{class}_\omega : \mathrm{Paths}(X) \to \mathscr{Y}.$$

DEFINITION (PREQUANTUM GROUPOID). *For a connected and simply connected diffeological space* $(X, \omega)$ *with discrete periods* $P_\omega$*, the* prequantum groupoid $\mathbf{T}_\omega$ *is the diffeological groupoid defined by:*

$$\mathrm{Obj}(\mathbf{T}_\omega) = X \quad \textit{and} \quad \mathrm{Mor}(\mathbf{T}_\omega) = \mathscr{Y} = \mathrm{Paths}(X)/\underset{\omega}{\sim}.$$

[24] For detailed technical exposition, consult [PIZ13, §8.42].

*The source and target maps* ends : $\mathscr{Y} \to \mathrm{X}\times\mathrm{X}$ *are induced by the endpoint maps on* Paths(X). *The groupoid composition is defined by the concatenation of paths:*

$$\mathrm{class}_\omega(\gamma)\cdot\mathrm{class}_\omega(\gamma') = \mathrm{class}_\omega(\gamma\vee\gamma'),$$

*for composable paths* $\gamma, \gamma'$. *The inverse is* $\mathrm{class}_\omega(\gamma)^{-1} = \mathrm{class}_\omega(\bar\gamma)$.

Let us introduce the operations of precomposition L and postcomposition R on $\mathscr{Y} = \mathrm{Mor}(\mathbf{T}_\omega)$:

$$\mathrm{L}(\tau)(\tau') = \mathrm{R}(\tau')(\tau) = \tau\cdot\tau'$$

where $\tau = \mathrm{class}_\omega(\gamma)$, $\tau' = \mathrm{class}_\omega(\gamma')$, and $\gamma(1) = \gamma'(0)$. The central property of this construction is articulated in the following theorem:

THEOREM (PREQUANTUM PRIMITIVE). *There exists a* 1*-form* $\boldsymbol{\lambda}$ *on* $\mathscr{Y}$ *satisfying:*

(1) $\boldsymbol{\lambda}$ *is the pushforward of* $\mathsf{K}\omega$, *i.e.*, $\mathrm{class}_\omega^*(\boldsymbol{\lambda}) = \mathsf{K}\omega$.
(2) $\boldsymbol{\lambda}$ *is invariant under precomposition and postcomposition in* $\mathscr{Y}$, *expressed as:*

$$\left\{\begin{array}{lcl} \mathrm{L}(\tau)^*(\boldsymbol{\lambda}\restriction \mathrm{Mor}_{\mathbf{T}_\omega}(x,\star)) & = & \boldsymbol{\lambda}\restriction \mathrm{Mor}_{\mathbf{T}_\omega}(x',\star), \\ \mathrm{R}(\tau)^*(\boldsymbol{\lambda}\restriction \mathrm{Mor}_{\mathbf{T}_\omega}(\star,x')) & = & \boldsymbol{\lambda}\restriction \mathrm{Mor}_{\mathbf{T}_\omega}(\star,x), \end{array}\right.$$

*We designate* $(\mathbf{T}_\omega, \boldsymbol{\lambda})$ *as the* prequantum groupoid *associated with* $(\mathrm{X},\omega)$, *and* $\boldsymbol{\lambda}$ *as the* prequantum primitive *of* $\omega$.

*The isotropy group at a point* $x$ *is equivalent to the torus of periods* $\mathrm{T}_\omega$*:*

$$\mathrm{Mor}_{\mathbf{T}_\omega}(x,x) = \mathrm{Loops}(\mathrm{X},x)/\underset{\omega}{\sim} \equiv \mathbf{R}/\mathrm{P}_\omega = \mathrm{T}_\omega.$$

*Sketch of Proof.* The proof, while technical, presents no significant difficulty [PIZ13, §8.42]. It is noteworthy that for the isotropy groups:

$$\int_{\hat x}^{\ell\vee\bar\ell'} \mathsf{K}\omega = \int_{\hat x}^{\ell} \mathsf{K}\omega - \int_{\hat x}^{\ell'} \mathsf{K}\omega \mod \mathrm{P}_\omega.$$

Consequently,

$$\ell \underset{\omega}{\sim} \ell' \iff \int_{\hat x}^{\ell} \mathsf{K}\omega = \int_{\hat x}^{\ell'} \mathsf{K}\omega \mod \mathrm{P}_\omega.$$

Moreover, the integral $\ell \mapsto \int_{\hat x}^{\ell} \mathsf{K}\omega$, from Loops(X, $x$) to $\mathbf{R}$, induces an isomorphism of diffeological groups from $\mathrm{Mor}_{\mathbf{T}_\omega}(x,x) = \mathrm{Loops}(\mathrm{X},x)/\underset{\omega}{\sim}$ to $\mathbf{R}/\mathrm{P}_\omega = \mathrm{T}_\omega$. ▶

DEFINITION (PREQUANTUM BUNDLE). *For each* $x \in \mathrm{X}$, *define* $\mathscr{Y}_x = \mathrm{Mor}_{\mathbf{T}_\omega}(x,\star)$ *equipped with the subset diffeology. The projection* $\hat\imath : \mathscr{Y}_x \to \mathrm{X}$ *is a principal bundle with structure group the isotropy group at* $x$, *equivalent to the torus of periods* $\mathrm{T}_\omega$. *The prequantum* 1*-form* $\boldsymbol{\lambda}$ *restricted to* $\mathscr{Y}_x$ *is a* connection 1-form $\lambda$ *whose curvature is* $\omega$, *i.e.*, $\hat\imath^*(\omega) = d\lambda$.

$$\begin{array}{ccccc} \mathrm{Paths}(\mathrm{X}) & \xrightarrow{\ \mathrm{class}_\omega\ } & \mathscr{Y} & \xleftarrow{\ \ j\ \ } & \mathscr{Y}_x \\ & \searrow^{\mathrm{ends}} & \downarrow{\scriptstyle \mathrm{ends}} & & \downarrow{\scriptstyle \hat\imath} \\ & & \mathrm{X}\times\mathrm{X} & \xrightarrow[pr_2]{\qquad} & \mathrm{X} \end{array}$$

Recall that a connection 1-form $\lambda$ on a principal bundle $\pi : Y \to X$ with structure group an irrational torus is a differential 1-form which is:

(1) Invariant under the torus action: $\tau^*(\lambda) = \lambda$ for all elements $\tau$ in the torus.
(2) Calibrated: $\hat{y}^*(\lambda) = \theta$, where $\hat{y}$ is the orbit map $\hat{y}(\tau) = \tau(y)$, and $\theta$ is the standard volume form on the torus.

Its exterior differential $d\lambda$ projects to the base space as a parasymplectic form, termed its *curvature*.

(1) For manifolds X, the construction of the prequantum bundle has been approached differently, including the non-simply connected case [PIZ95]. The comprehensive result is that such generalized prequantum bundles are classified by the extension group $\mathrm{Ext}(\mathbf{H}^1(X, \mathbf{Z}), P_\omega)$. The general non-simply connected case within diffeology is an ongoing area of research.
(2) The prequantum groupoid, constructed as a quotient of $(\mathrm{Paths}(X), \mathbf{K}\omega)$ and endowed with its prequantum 1-form $\boldsymbol{\lambda}$, exhibits a source-target symmetry absent in the prequantum bundle. This symmetry suggests a potential deep connection between Feynman path integral quantization and geometric quantization, a relationship warranting further investigation.
(3) In this general picture a *wave function* would be a morphism $\psi$ from the prequantum groupoid $\mathscr{Y}$, with values into some Abelian algebraic structure, that is, $\psi(\tau \cdot \tau') = \psi(\tau) \cdot \psi(\tau')$, where $\tau = \mathrm{class}_\omega(\gamma)$ and $\tau' = \mathrm{class}_\omega(\gamma')$ with $\gamma(1) = \gamma'(0)$, and denoted multiplicatively for coherence with the familiar case when $P_\omega = h\mathbf{Z}$, and the algebraic structure is $\mathbf{C}$.
(4) Focusing solely on the prequantum bundle limits the study of symmetries of $\omega$. While a Lie group G preserving $\omega$ can be represented in the group of quantomorphisms only under specific conditions (see, for example, [Sou70]), it admits a direct representation in $\mathrm{Aut}(\mathscr{Y}, \boldsymbol{\lambda})$ via $\underline{g}(\mathrm{class}_\omega(\gamma)) = \mathrm{class}_\omega(\underline{g} \circ \gamma)$. The moment map associated with this action is the projection of the paths moment map $\Psi : \mathrm{class}_\omega(\gamma) \to \mathscr{G}^*$, which then descends to $\psi : X \times X \to \mathscr{G}^*/\Gamma$.

*In summary*, we have constructed a prequantum groupoid $(\mathbf{T}_\omega, \boldsymbol{\lambda})$ for connected and simply connected diffeological spaces $(X, \omega)$ with discrete periods. This object, built intrinsically from the geometry of the classical system via a diffeological quotient of the space of paths in X using the chain-homotopy operator, serves as a unified structure for prequantization. The prequantum groupoid $(\mathbf{T}_\omega, \boldsymbol{\lambda})$ generalizes the traditional prequantum bundle, which arises as a specific projection $\mathscr{Y}_x = \mathrm{Mor}_{\mathbf{T}_\omega}(x, \star) \to X$ for a fixed base point $x$. The groupoid structure reveals a unique source-target symmetry absent in the prequantum bundle, hinting at a deeper connection between Feynman path integral quantization and geometric quantization, a relationship deserving of further exploration.

Furthermore, we have demonstrated that the prequantum groupoid provides a more comprehensive approach to understanding symmetries of the symplectic form $\omega$, overcoming limitations inherent in the prequantum bundle formulation. The entire symmetry group $\mathrm{Diff}(X, \omega)$ admits a direct representation as automorphisms of $(\mathbf{T}_\omega, \boldsymbol{\lambda})$ via $\underline{g}(\text{class}_\omega(\gamma)) = \text{class}_\omega(\underline{g}\circ\gamma)$, and the associated paths moment map offers powerful tools for exploiting these symmetries.

This generalization not only addresses the need to extend geometric quantization to spaces beyond smooth manifolds but also opens new avenues for investigating the fundamental connections between classical and quantum mechanics. The prequantum groupoid, with its enhanced symmetry properties and broader applicability, provides a robust foundation for future research in quantization theory, particularly in contexts involving infinite-dimensional spaces, as in quantum field theory, and spaces with singularities, as in symplectic reduction. The ongoing work on the non-simply connected case promises to further solidify the relevance and efficacy of this generalized approach.

**23. Every Symplectic Manifold is a Coadjoint Orbit.** In the late 1960s, Kostant, Kirillov, and Souriau, working from different perspectives, demonstrated that a symplectic manifold $(M,\omega)$ homogeneous under the action of a Lie group is isomorphic to a coadjoint orbit [Kos70, Sou70, Kir74].[25] Souriau's proof relied on the moment map, a fundamental tool he introduced during that period.

Consider a connected symplectic manifold $(M,\omega)$. Its group of automorphisms, denoted $\text{Diff}(M,\omega)$ (also known as the group of symplectomorphisms), acts transitively on $M$. It is natural to consider $(M,\omega)$ as a coadjoint orbit of its group of symplectomorphisms, analogous to the Kostant-Kirillov-Souriau (KKS) theorem. However, $\text{Diff}(M,\omega)$ is not a Lie group, and the KKS theorem cannot be directly applied. This issue has not been extensively addressed, perhaps due to the lack of a suitable framework. Functional analysis could be considered, but its heavy formalism is less suited to the intuitive geometric nature of this problem.

As we will see in this section, diffeology provides a natural and efficient framework to address this, achieving a significant result with minimal additional structure. Thanks to the universal moment map discussed earlier, $(M,\omega)$ is equivalent to an affine coadjoint orbit of $\text{Diff}(M,\omega)$. Although this affine action, defined by a twisted cocycle of the symplectomorphisms, is a limitation, we prefer to identify the symplectic manifold with an ordinary coadjoint orbit, i.e., an orbit of the standard linear coadjoint action. This can be achieved by passing to a central extension of the group of automorphisms, using the prequantum bundle construction.

It is worth noting that students of symplectic geometry often find the fundamental nature of symplectic manifolds confusing. I have even encountered a question on MathOverflow asking whether every symplectic manifold is a cotangent bundle. With the construction enabled by the flexible diffeological framework, we can provide a clear answer: No, *every symplectic manifold is a coadjoint orbit.*

Instead of the entire group of symplectomorphisms, we consider the subgroup of *Hamiltonian diffeomorphisms*, which is also transitive on $(M,\omega)$ [Boo69]. It is useful to introduce this group in the context of any parasymplectic space, as its definition does not depend on the space being a manifold.

[25] Or a covering of a coadjoint orbit.

THEOREM (HAMILTONIAN DIFFEOMORPHISMS). *For any connected parasymplectic space* $(X,\omega)$*, there exists a* largest connected subgroup $\mathrm{Ham}(X,\omega) \subset \mathrm{Diff}(X,\omega)$ *whose holonomy is trivial, meaning that the values of its paths moment map on loops are zero. This subgroup is called the group of* Hamiltonian diffeomorphisms.

In what follows, we denote $\mathrm{Diff}(X,\omega)$ by $G_\omega$ and $\mathrm{Ham}(X,\omega)$ by $H_\omega$.

*Sketch of Proof.* Let $\tilde{G}^0_\omega$ be the universal covering of the identity component $G^0_\omega$ of $G_\omega$, let $\pi : \tilde{G}^0_\omega \to G^0_\omega$ be the projection. Every element $\gamma$ of the holonomy group $\Gamma_\omega$ is a closed 1-form on $G_\omega$. Let $k(\gamma)$ be the real homomorphism on $\tilde{G}^0_\omega$ such that $\pi^*(\gamma) = d[k(\gamma)]$. Let

$$\hat{H}_\omega = \bigcap_{\gamma \in \Gamma_\omega} \ker\big(k(\gamma)\big).$$

Then, $H_\omega = \pi(\hat{H}^0_\omega)$ is the largest connected subgroup of $G_\omega$ with trivial holonomy. We denote by $\pi : \tilde{G}^0_\omega \to G^0_\omega$ the projection from the universal covering. ►

This theorem and its proof rely on several features of diffeology:

* The existence of a universal simply connected covering for all diffeological spaces.
* The consistent definition of differential forms on all diffeological spaces.
* The consistent definitions of pullback and exterior differential.
* The property that a closed 1-form on a simply connected diffeological space is exact.

Even when X is a manifold, proving this theorem in the traditional framework requires the ad hoc construction of time-dependent Hamiltonians [PIZ13, §9.16].

Next, let $\bar{\mu}_\omega$ denote the moment map for the group $H_\omega$. Then, based on the preceding discussion of moment maps, we have:

$$\bar{\mu}_\omega : X \to \mathscr{H}^*, \quad \text{and} \quad \bar{\mu}_\omega(g(x)) = \mathrm{Ad}_*(g)(\bar{\mu}_\omega(x)) + \bar{\theta}_\omega(g),$$

where $\bar{\theta}_\omega$ is the cocycle associated with the action of $H_\omega$.

The fundamental theorem upon which the result that every symplectic manifold is a coadjoint orbit rests is the following [PIZ17]:

THEOREM (TO BE A SYMPLECTIC FORM). *Let* $\omega$ *be a parasymplectic form on a manifold* M. *Then,* $\omega$ *is symplectic if and only if:*

(1) *The action of* $H_\omega$ *is transitive on* M.
(2) *The associated moment map* $\bar{\mu}_\omega$ *is injective.*

*Sketch of Proof.* The condition that the action of $H_\omega$ is transitive is due to Boothby [Boo69]. The demonstration that $\bar{\mu}_\omega$ is injective is based on the following: let $m_0$ and $m_1$ be two distinct points of M. Let $p$ be a path connecting $m_0$ to $m_1$, thus $\bar{\mu}_\omega(m_1) - \bar{\mu}_\omega(m_0) = \bar{\Psi}_\omega(p)$. Since M is Hausdorff, there exists a smooth real function $f \in \mathscr{C}^\infty(M,\mathbf{R})$ with compact support such that $f(m_0) = 0$ and $f(m_1) = 1$. Let $\xi$ denote the symplectic gradient field associated with $f$, and let F be the exponential of $\xi$. Thus, $\bar{\Psi}_\omega(p)(F) = dt$ by the application of the paths moment map to manifolds. Hence, $(\bar{\mu}_\omega(m_1) - \bar{\mu}_\omega(m_0))(F) = dt \neq 0$, implying $\bar{\mu}_\omega(m_0) \neq \bar{\mu}_\omega(m_1)$. Therefore, $\bar{\mu}_\omega$ is injective. ►

It is important to note that these two conditions are necessary. For example, on $\mathbf{R}^2$, the 2-form $\omega = (x^2 + y^2)dx \wedge dy$ is closed, and $\bar{\mu}_\omega$ is injective, but the action of $H_\omega$ is not transitive, since $\ker(\omega)_0 = \mathbf{R}^2$. Now, the orbit map $\hat{m} : H_\omega \to M$ is a subduction [Don84], i.e., the diffeology of M is the pushforward of the functional diffeology on $H_\omega$. Then:

THEOREM (SYMPLECTIC MANIFOLDS ARE COADJOINT ORBITS). *Let* $(M,\omega)$ *be a symplectic manifold. Then, the moment map* $\bar{\mu}_\omega : M \to \mathscr{H}^*$ *is a diffeomorphism from* $(M,\omega)$ *to its image* $\mathscr{O} \subset \mathscr{H}^*_\omega$*. Here,* $\mathscr{O}$ *is equipped with the quotient diffeology induced by the action of* $H_\omega$ *equipped with the functional diffeology, and its symplectic structure is the pushforward of* $\omega$*.*

However, the orbit $\mathscr{O}$ can be affine when the cocycle $\bar{\theta}_\omega$ is non-trivial. To eliminate this affine character, we consider a prequantum bundle $\pi : Y \to M$, with connection form $\lambda$ and structure group $T_\omega$, where $T_\omega = \mathbf{R}/P_\omega$ and $P_\omega$ is the group of periods of $\omega$. For manifolds, this construction applies to any homotopy type of M [PIZ95].

THEOREM (PREQUANTUM BUNDLE AUTOMORPHISMS). *The identity component of the group of automorphisms of the prequantum bundle* $(Y,\lambda)$ *is a central extension of the group of Hamiltonian diffeomorphisms of* $(M,\omega)$ *by the torus of periods* $T_\omega$*, summarized by the following exact sequence:*

$$\mathbf{1} \longrightarrow T_\omega \longrightarrow \mathrm{Aut}(Y,\lambda)^0 \xrightarrow{\mathrm{pr}} \mathrm{Ham}(M,\omega) \longrightarrow \mathbf{1}$$

The projection pr induces a pullback $\mathrm{pr}^* : \mathscr{H}^*_\omega \to \mathscr{A}^*_\lambda$, where $\mathscr{A}^*_\lambda$ denotes the space of momenta of $\mathrm{Aut}(Y,\lambda)$. Let $\mu_Y : Y \to \mathscr{A}^*_\lambda$ be the moment map associated with $d\lambda = \pi^*(\omega)$, i.e., $\mu_Y(y) = \hat{y}^*(\lambda)$. Note that the moment map $\mu_Y$ is constant on the fibers of $\pi : Y \to M$, since $\mathrm{Aut}(Y,\lambda)$ is a central extension by $T_\omega$. Let $\mu_M$ be its projection onto M, i.e., $\mu_Y = \mu_M \circ \pi$. Then:

THEOREM (SYMPLECTIC MANIFOLDS AS LINEAR COADJOINT ORBITS). *The moment map* $\mu_M$*, of the action of the group of automorphisms of* $(Y,\lambda)$ *on* M*, is injective and identifies* M *with a linear (i.e., non-affine)* $\mathrm{Ad}_*$*-coadjoint orbit.*

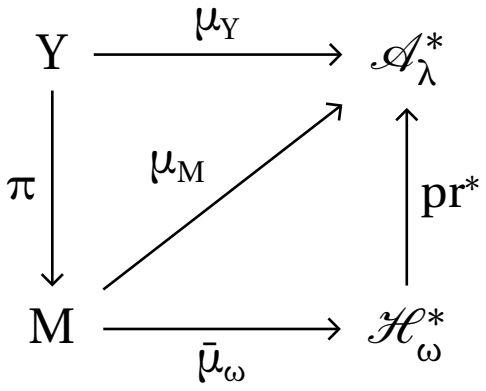

*In conclusion*: to obtain symplectic manifolds as *linear* coadjoint orbits, we leverage the capabilities of diffeology to handle infinite-dimensional groups (the group of quantomorphisms) and singular spaces (the torus of periods). This example demonstrates how these combined capabilities of diffeology provide a comprehensive and satisfactory resolution to a problem difficult to address in traditional frameworks.

**24. Application: A Foundational Reformulation of Prequantization.** Geometric quantization, the bridge between classical and quantum mechanics, has long been challenged by singular and infinite-dimensional spaces, and by the mysterious appearance of central extensions of symmetry groups. In [PIZ25b], we use diffeology to reformulate the very foundations of prequantization, resolving these issues with a more fundamental object.

Instead of positing a prequantum bundle, we construct a **prequantum groupoid** directly from the path-space geometry of the classical system $(X,\omega)$. This path-based approach, resonant with Feynman's philosophy, yields two profound insights. First, the "quantum phase" group (the circle $T_\omega$) is not an external object put in by hand; it **emerges naturally** as the isotropy group of the groupoid, constructed from the geometry of loops in the space.

Second, we prove that the full symmetry group $\mathrm{Diff}(X, \omega)$ acts **isomorphically** as automorphisms of the entire prequantum groupoid structure. This demonstrates that the famous central extensions are not a fundamental feature of quantization itself, but rather an artifact that appears when one descends from this more complete groupoid structure to a less fundamental object, such as a bundle or a space of sections over X. This framework provides a conceptually clearer, more general, and more powerful starting point for quantization, applicable to the very spaces where physics needs it most.

## 7. Diffeology Meets Noncommutative Geometry

The connection between diffeology and Alain Connes' Noncommutative Geometry has a rich history, dating back to the early 1980s. Connes introduced noncommutative geometry [AC80, AC95] to address the spectrum of the Hamiltonian in the quantum model of quasi-periodic potentials [DS75].

During the same period, Jean-Marie Souriau laid the groundwork for diffeological spaces with his "groupes différentiels" [Sou80], initially developed specifically for studying groups of diffeomorphisms. It wasn't until 1983 that, intrigued by Connes' work on irrational foliations, with Paul Donato, we discovered how Souriau's "groupes différentiels" provided a surprisingly non-trivial description of the topologically trivial irrational torus [DIZ83]. This pivotal moment marked the shift from Souriau's "groupes différentiels" to the theory of diffeological spaces and spurred its development as a generalization of differential geometry. This included a significant development of diffeological groups and Cartan-de Rham calculus, as well as the introduction of new areas of study such as fiber bundle theory, higher homotopy theory, and the study of singular objects like orbifolds, quasifolds, and spaces with internal singularities (e.g., Klein stratifications).

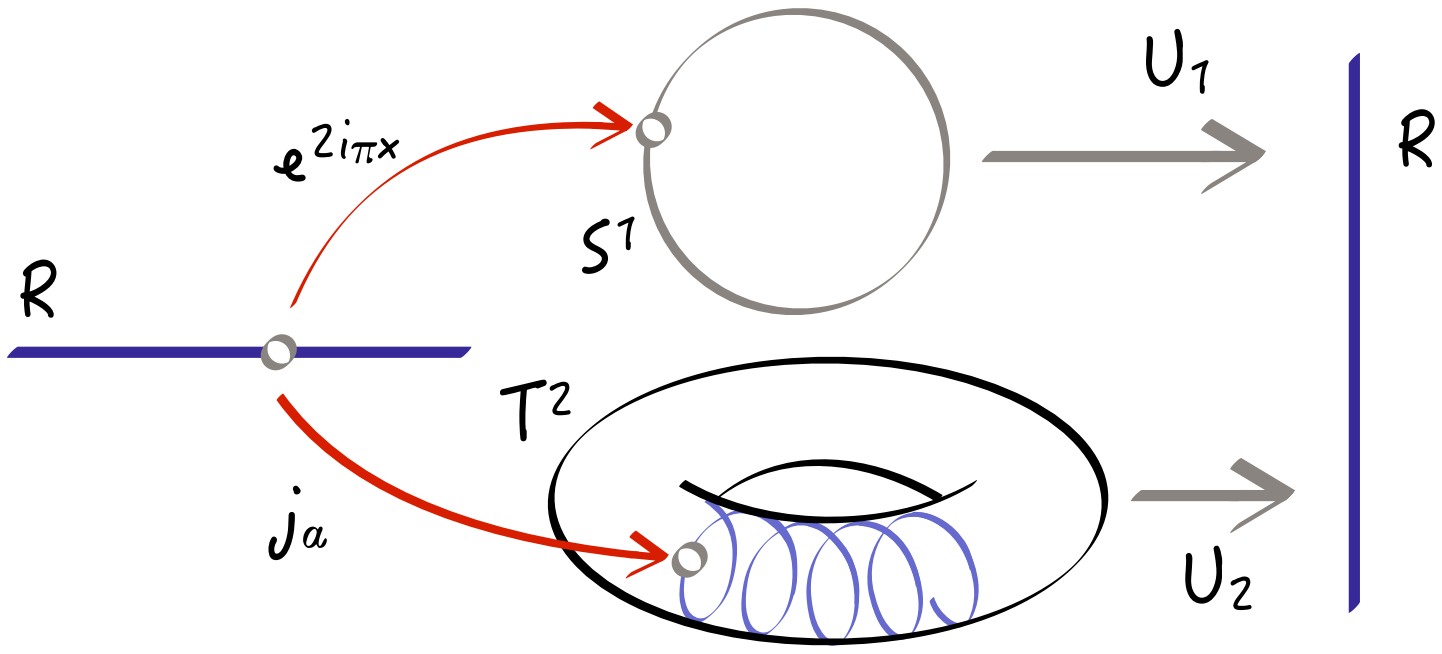

Figure 3. Periodic and quasi-periodic potentials.

However, the relationship between diffeology and noncommutative geometry remained unclear for many years. It is only recently that the deep connection between these two fields has been revealed through the construction of $\mathbf{C}^*$-algebras structurally associated with orbifolds and quasifolds [IZL18, IZP21]. This construction has the minimal property that diffeomorphic quasifolds correspond to Morita equivalent C*-algebras, unveiling the organic relationship between diffeology and noncommutative geometry.

The connection between diffeology and noncommutative geometry is established following Jean Renault's groupoid approach to noncommutative geometry. It concerns exclusively orbifolds and quasifolds, as defined in §9. It follows the scheme:

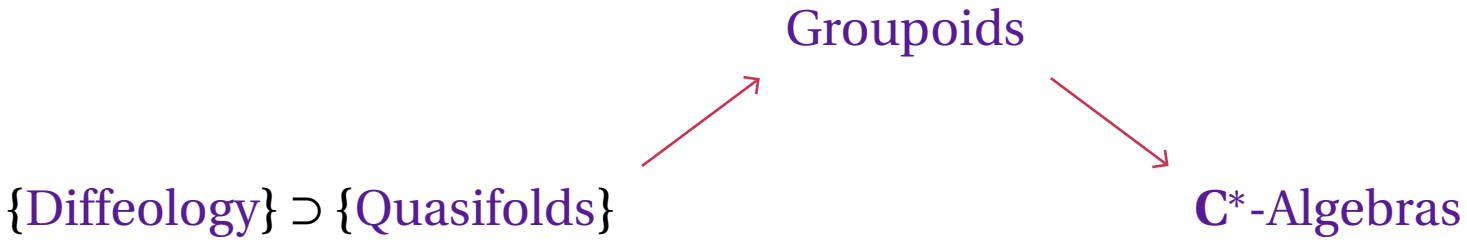

Let $\mathcal{Q}$ be a quasifold, which includes orbifolds as a particular case. The structure groupoid $\mathbf{G}$, associated with an atlas $\mathcal{A}$, with a strict generating family $\mathcal{F}$ and nebula $\mathcal{N}$ as the object space, satisfies Renault's conditions:

DEFINITION ($\mathbf{C}^*$-ALGEBRA OF A QUASIFOLD). *The structure groupoid of a quasifold $\mathcal{Q}$, as defined in §9, is étale and Hausdorff. As such, it meets the conditions required for Renault's construction* [Ren80] *of a $\mathbf{C}^*$-algebra $\mathfrak{A}$, obtained by completing the compactly supported complex functions on $\mathbf{G}$ with respect to the uniform norm, equipped with the counting measure. The convolution and involution are defined by:*

$$f * g(\gamma) = \sum_{\beta \in \mathbf{G}^x} f(\beta \cdot \gamma) g(\beta^{-1}) \quad \textit{and} \quad f^*(\gamma) = f(\gamma^{-1})^*,$$

*where $x = \mathrm{src}(\gamma)$, and $z^*$ is the complex conjugate of $z \in \mathbf{C}$. The sums are assumed to converge.*

The association of the $\mathbf{C}^*$-algebra $\mathfrak{A}$ to the atlas $\mathcal{A}$ is a morphism in the sense:

THEOREM (EQUIVALENT $\mathbf{C}^*$-ALGEBRAS). *Different atlases yield Morita-equivalent algebras, and $\mathbf{C}^*$-algebras associated with diffeomorphic quasifolds are Morita equivalent.*

*Sketch of Proof.* The proof utilizes the Muhly-Renault-Williams theorem on $\mathbf{C}^*$-algebras defined by groupoids [MRW87], by constructing a form of "cobordism" between the two groupoids, composed of arrows of the same kind. The elements of this "cobordism" are illustrated in Figure (4). They are germs of equivariant transversal diffeomorphisms $\varphi$ between atlases, lifting local diffeomorphisms $\psi$ over a germ of the identity of $\mathcal{Q}$. ►

The Morita equivalence between $\mathbf{C}^*$-algebras is the minimal condition for a categorical functor from the subcategory of quasifolds to the category of $\mathbf{C}^*$-algebras, as required by noncommutative geometry.

*Examples.* Let us conclude with two examples, apart from the well-known irrational torus associated with the irrational rotation algebra [Rie81]. The first example is the $*$-algebra of the orbifold $\Delta_1 = \mathbf{R}/\{\pm 1\}$. The singleton $\mathcal{F} = \{\mathrm{sq} : t \mapsto t^2\}$ is a strict generating family, and the structure groupoid $\mathbf{G}$ is the groupoid of the action of $\Gamma = \{\pm 1\}$, i.e.,

$$\mathrm{Obj}(\mathbf{G}) = \mathbf{R} \ \text{ and } \ \mathrm{Mor}(\mathbf{G}) = \{(t, \varepsilon, \varepsilon t) \mid \varepsilon = \pm 1\} \simeq \mathbf{R} \times \{\pm 1\}.$$

A continuous function $f$ on $\mathrm{Mor}(\mathbf{G})$ to $\mathbf{C}$ is a pair of functions $f = (a, b)$, where $a(t) = f(t, 1)$ and $b(t) = f(t, -1)$. With this convention, applying the definition reminded in the preamble, the convolution product is represented by matrix multiplication in the module $M_2(\mathbf{C}) \otimes \mathscr{C}^0(\mathbf{R}, \mathbf{C})$, with

$$f = (a, b) \mapsto \mathrm{M} \quad \text{with} \quad \mathrm{M}(t) = \begin{pmatrix} a(t) & b(-t) \\ b(t) & a(-t) \end{pmatrix} \quad \text{and} \quad \mathrm{M}^*(t) = [{}^\tau \mathrm{M}(t)]^*,$$

where the superscript $\tau$ denotes transposition, and the asterisk denotes element-wise complex conjugation.

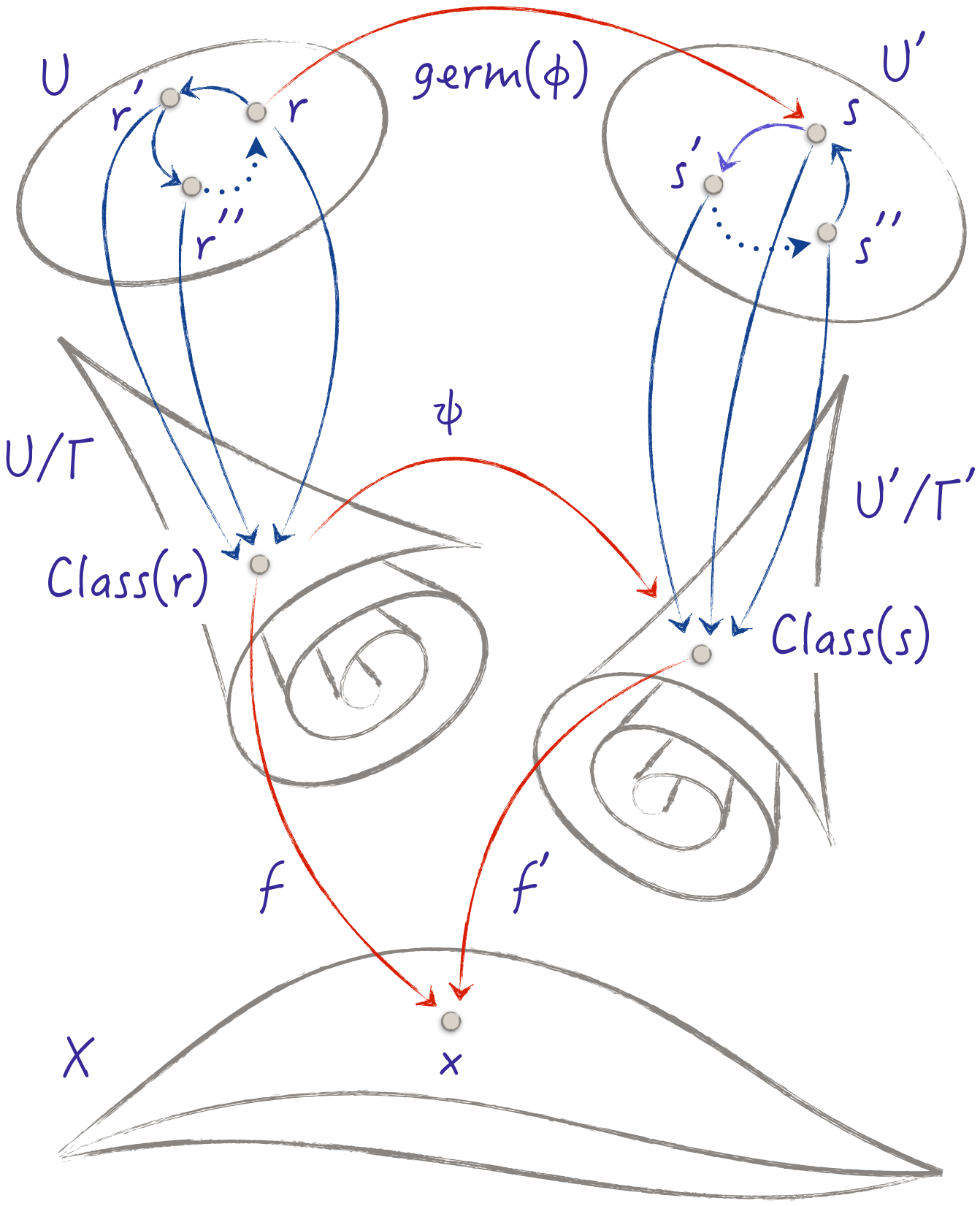

Figure 4. The MRW "Cobordism".

For the second example, consider the quasifold $\mathscr{Q} = \mathbf{R}/\mathbf{Q}$. The structure groupoid is the groupoid of the action of $\mathbf{Q}$ on $\mathbf{R}$. Its arrow space is simply $\mathbf{R} \times \mathbf{Q}$. Now, the quotient $\mathbf{R}/\mathbf{Q}$ is also diffeomorphic to the $\mathbf{Q}$-circle $S_{\mathbf{Q}} = S^1/\mathscr{U}_{\mathbf{Q}}$, where $\mathscr{U}_{\mathbf{Q}} = \{\exp(2i\pi r)\}_{r \in \mathbf{Q}}$ is the subgroup of rational roots of unity. The arrow space of the structure groupoid of this quotient is simply $S^1 \times \mathscr{U}_{\mathbf{Q}}$. The exponential $x \mapsto z = e^{2i\pi x}$ establishes a MRW-equivalence between the two groupoids. Their associated algebras are Morita equivalent. The associated algebra, which we will denote by $\mathfrak{S}$, consists of families of continuous complex functions indexed by rational roots of unity, i.e., $(f_\tau)_{\tau \in \mathscr{U}_{\mathbf{Q}}}$, with finite support. The convolution and algebra conjugation are:

$$(f * g)_\tau(z) = \sum_\sigma f_{\bar{\sigma}\tau}(\sigma z) g_\sigma(z) \quad \text{and} \quad f^*_\tau(z) = f_{\bar{\tau}}(\tau z)^*,$$

where $\bar\tau = 1/\tau = \tau^*$, the complex conjugate.

Now, consider $f$ and let $\mathscr{U}_p$ be the subgroup in $\mathscr{U}_{\mathbf{Q}}$ generated by $\mathrm{Supp}(f)$; this is the group of some root of unity $\varepsilon$ of order $p$. Let $\mathrm{M}_p(\mathbf{C})$ be the space of $p \times p$ matrices with complex coefficients. Define $f \mapsto \mathrm{M}$, from $\mathfrak{S}$ to $\mathrm{M}_p(\mathbf{C}) \otimes \mathscr{C}^0(\mathrm{S}^1, \mathbf{C})$, by

$$\mathrm{M}(z)^\sigma_\tau = f_{\bar\sigma\tau}(\sigma z), \quad \text{for all } z \in \mathrm{S}^1 \text{ and } \sigma, \tau \in \mathscr{U}_p.$$

This gives a representation of $\mathfrak{S}$ in the tensor product of the space of endomorphisms of the infinite-dimensional $\mathbf{C}$-vector space $\mathrm{Maps}(\mathscr{U}_{\mathbf{Q}}, \mathbf{C})$ by $\mathscr{C}^0(\mathrm{S}^1, \mathbf{C})$, with finite support.

*In conclusion*, A key distinction between diffeology and noncommutative geometry lies in their respective approaches to geometric objects. Noncommutative geometry frequently prioritizes algebraic constructions, sometimes at the expense of a clear underlying geometric picture. In contrast, diffeology, with its emphasis on the construction of diffeological spaces, offers a more geometrically driven methodology. For example, instead of beginning with the irrational rotation algebra $\mathrm{A}_\theta$ as a purely algebraic entity, diffeology first constructs the irrational torus $\mathrm{T}_\theta$ as a well-defined diffeological space, and subsequently derives the associated $\mathbf{C}^*$-algebra through a categorical process. This approach provides a clearer geometric foundation, aligning with the perspective of a geometer.

## 8. Conclusion: The Power and Potential of Diffeology

From the inscription above Plato's Lyceum, 'Let no one ignorant of geometry enter here', to the modern era, physics has been inextricably linked with geometry. The evolution of physics has often mirrored the development of geometric tools, with differential geometry playing a pivotal role. Yet, classical differential geometry encounters limitations when faced with complex spaces naturally arising in modern mathematics and physics, such as quotients, infinite-dimensional spaces, and spaces with singularities.

Diffeology emerges as a natural and powerful extension, offering a simple yet profoundly versatile framework. By defining smoothness intrinsically via plots, it allows us to treat singular spaces as autonomous geometrical objects, endowed with their own rich structure. This aligns perfectly with Felix Klein's vision of geometry: the study of a space through the lens of its transformation group and the corresponding invariants. For a diffeological space X, this group is its group of diffeomorphisms Diff(X), and invariants like dimension, Klein strata, homotopy groups, cohomology groups, and classifying groups like $\mathbf{Fl}(\mathrm{X}, \mathbf{R})$ reveal its underlying geometry, see [PIZ25a].

As demonstrated throughout this paper — from the arithmetic encoded in the diffeology of the irrational torus, to the rigorous handling of orbifolds and the definitive geometry of orbit spaces — diffeology provides not just solutions, but often more direct and conceptually clearer paths than traditional methods constrained by manifold assumptions. It offers a unifying language capable of expressing complex geometric and physical ideas with fidelity and rigor.

Yet, the ultimate measure of a theory's vitality lies not only in the old problems it solves, but in the new and fruitful questions it generates. Each solution presented here opens a new avenue for research. The 'orthofold' is not merely an answer but a new class of geometric objects demanding its own theory of invariants and structures. The 'prequantum groupoid' is not just a technical fix but a fundamental object whose representation theory remains to be explored. Finally, the successful treatment of singular symplectic reduction in our examples is not an end-point; it is a proof of concept that demands the formulation of a complete, general theory of reduction for all diffeological spaces.

This capacity to transform foundational puzzles into fertile new research programs stems from a commitment to a certain mathematical *aisthēsis* [26]: the pursuit of the simplest, most natural structures that reveal the essence of a problem, rather than a display of technical virtuosity for its own sake. This approach — providing not just a collection of tools, but a coherent framework that reveals the next set of questions — is perhaps the most compelling answer to the question: "Why Diffeology?".

## A Personal Note

I first met Vladimir Arnold forty-one years ago in Moscow, sitting in the kitchen of his apartment in Yugo-Zapadnaya. At the time, access to the university was restricted, and I could not obtain the necessary *propusk* to pass the gate. I had come to Moscow as a postdoc but at the Steklov Intitute, where I had finished my work on the classification of SO(3)-symplectic manifolds which I wanted to present to him. I remember vividly that he spoke to me for a long time—not about my work, but about his own passion at the moment: the isochrone pendulum. The journey from Ulitsa Gubkina to Yugo-Zapadnaya remains etched in my memory.

We met several times afterwards in France, developing a warm connection. He eventually came to regard me as his student, a title I accept with pride; he was, alongside Jean-Marie Souriau, one of my two great masters. It is a profound honor to publish this paper today in the journal that bears his name.

[26]From the Ancient Greek αἴσθησις (aísthēsis): the sensory experience, perception, and in a philosophical sense, the faculty of judgment and appreciation of what is fitting or true.

EINSTEIN INSTITUTE OF MATHEMATICS, THE HEBREW UNIVERSITY OF JERUSALEM, EDMOND J. SAFRA CAMPUS ON GIVAT RAM, 9190401 ISRAEL.

*Email address*: `piz@math.huji.ac.il`

*URL*: `http://math.huji.ac.il/~piz`